\documentclass[11pt]{amsart}
\usepackage{palatino}
\usepackage{amsfonts}
\usepackage{amssymb}
\usepackage{amscd}
\usepackage[dvips]{graphicx}

\theoremstyle{plain}

\newtheorem{theorem}{Theorem}[section]
\newtheorem{proposition}[theorem]{Proposition}
\newtheorem{corollary}[theorem]{Corollary}

\theoremstyle{definition}
\newtheorem{definition}[theorem]{Definition}
\newtheorem{remark}[theorem]{Remark}

\newtheorem{conjecture/question}[theorem]{Conjecture/Question}

\newtheorem{remark/definition}[theorem]{Remark/Definition}
\newtheorem{terminology/notation}[theorem]{Terminology/Notation}

\newcommand{\marginlabel}[1]%
  {\mbox{}\marginpar{\raggedleft\hspace{0pt}\bfseries\sf#1}}

\def\PP{{\textbf P}}
\def\OO{\mathcal{O}}
\def\cN{\mathcal{N}}

\def\cB{\mathcal{B}}
\def\cA{\mathcal{A}}
\def\F{\mathcal{F}}
\def\P{\mathcal{P}}

\def\L{\mathcal{L}}

\def\cM{\mathcal{M}}
\def\cR{\mathcal{R}}

\def\cU{\mathcal{U}}
\def\cC{\mathcal{C}}

\def\Pic0{{\rm Pic}^0(X)}

\def\mm{\overline{\mathcal{M}}}

\def\gp{\overline{\mathcal{GP}}}
\def\cgp{\mathcal{GP}}

\theoremstyle{remark}

\begin{document}

\title{Rational maps between moduli spaces of curves and Gieseker-Petri divisors }
\author[G. Farkas]{Gavril Farkas}
\address{Humboldt-Universit\"at zu Berlin, Institut f\"ur
Mathematik,
 10099 Berlin} \email{{\tt farkas@math.hu-berlin.de}}

\thanks{Research  partially supported by an Alfred P. Sloan Fellowship, the NSF Grants DMS-0450670 and DMS-0500747
and a 2006 Texas Summer Research Assignment. Most of this paper has
been written while visiting the Institut Mittag-Leffler in
Djursholm in the Spring of 2007. Support from the institute is gratefully acknowledged.}

\maketitle

For a a general smooth projective curve $[C]\in \cM_g$ and an arbitrary
line bundle $L\in \mbox{Pic}(C)$, the
\emph{Gieseker-Petri} theorem states that the multiplication map
$$\mu_0(L):H^0(C, L)\otimes H^0(C, K_C\otimes L^{\vee})\rightarrow
H^0(C, K_C)$$ is injective. The
theorem, conjectured by Petri and proved by Gieseker
\cite{G}  (see \cite{EH3} for a much simplified proof), lies at the
cornerstone of the theory of algebraic curves. It implies
that the variety $G^r_d(C)=\{(L, V): L\in \mbox{Pic}^d(C),
V\in G(r+1, H^0(L))\}$ of linear series of degree $d$ and dimension
$r$ is smooth and of expected dimension $\rho(g, r, d):=g-(r+1)(g-d+r)$ and
that the forgetful map $G^r_d(C)\rightarrow W^r_d(C)$ is a rational
resolution of singularities (see \cite{ACGH} for many other
applications). It is an old open problem to describe the locus
$\mathcal{GP}_g\subset \cM_g$ consisting of curves $[C]\in \cM_g$
such that there exists a line bundle $L$ on $C$ for which the
Gieseker-Petri theorem fails. Obviously $\mathcal{GP}_g$ breaks up
into irreducible components depending on the numerical types of
linear series. For fixed integers $d, r\geq 1$ such that $g-d+r\geq 2$, we define the locus
$\mathcal{GP}_{g, d}^r$ consisting of curves $[C]\in \cM_g$ such
that there exist a pair of linear series $(L, V)\in G^r_d(C)$ and
$(K_C\otimes L^{\vee}, W)\in G_{2g-2-d}^{g-d+r-1}(C)$ for which the
multiplication map
$$\mu_0(V, W): V\otimes W\rightarrow H^0(C, K_C)$$ is not injective. Even though certain components of $\mathcal{GP}_g$ are
well-understood, its global geometry seems exceedingly
complicated. If $\rho(g, r, d)=-1$, then $\mathcal{GP}_{g, d}^r$
coincides with the Brill-Noether divisor $\cM_{g, d}^r$ of curves
$[C]\in \cM_g$ with $G^r_d(C)\neq \emptyset$ which has been
studied by Eisenbud and Harris in \cite{EH2} and used to
prove that $\mm_g$ is of general type for $g\geq 24$. The locus
$\mathcal{GP}_{g, g-1}^1$ can be identified with the divisor of
curves carrying a vanishing theta-null and this has been studied by
Teixidor (cf. \cite{T}). We proved in \cite{F2} that
for $r=1$ and $(g+2)/2\leq d\leq g-1$, the locus $\mathcal{GP}_{g,
d}^1$ always carries a divisorial component. It is conjectured  that
the locus $\mathcal{GP}_g$ is pure of codimension $1$ in $\cM_g$ and
we go some way towards proving this conjecture. Precisely, we show that
$\mathcal{GP}_g$ is supported in codimension $1$ for every possible
numerical type of a linear series:
\begin{theorem}\label{gisp1}
For any positive integers $g, d$ and $r$ such that $\rho(g, r,
d)\geq 0$ and $g-d+r\geq 2$, the locus $\mathcal{GP}_{g, d}^r$
 has a divisorial component in $\cM_g$.
\end{theorem}

The main issue we address in this paper
is a detailed intersection theoretic study of a rational map between two
different moduli spaces of curves. We fix  $g:=2s+1\geq 3$. Since
$\rho(2s+1, 1, s+2)=1$  we can define a rational map between moduli spaces of curves
$$\phi: \mm_{2s+1} -->\mm_{1+\frac{s}{s+1}{2s+2\choose s}}, \mbox{  } \mbox{ } \  \phi([C]):=[W^1_{s+2}(C)].$$
The fact that $\phi$ is well-defined, as well as a justification for
the formula of the genus  $g':=g(W^1_{s+2}(C))$ of the curve of
special divisors of type $\mathfrak g^1_{s+2}$, is given in Section
3. It is known that $\phi$ is generically injective (cf. \cite{PT},
\cite{CHT}). Since $\phi$ is the only-known rational map between two
moduli spaces of curves and one of the very few natural examples of
a rational map admitted by $\mm_g$, its study is clearly of
independent interest. In this paper we carry out a detailed
enumerative study of $\phi$ and among other things, we determine the
pull-back map $\phi^*: \mbox{Pic}(\mm_{g'})\rightarrow
\mbox{Pic}(\mm_g)$ (see Theorem \ref{divisorspullback} for a precise
statement). In particular we have the following formula concerning
slopes of divisor classes pulled back from $\mm_{g'}$ (For the
definition of the slope function $s:\mathrm{Eff}(\mm_g)\rightarrow
\mathbb R\cup \{\infty\}$ on the cone of effective divisors we refer
to \cite{HMo} or \cite{FP}):
\begin{theorem}\label{slopepull}
We set $g:=2s+1$ and $g':=1+\frac{s}{s+1}{2s+2\choose s}$. For any
divisor class $D\in \mathrm{Pic}(\mm_{g'})$ having slope $s(D)=c$,
we have the following formula for the slope of $\phi^*(D)\in
\mathrm{Pic}(\mm_g)$:
$$s(\phi^*(D))=6+ \frac{8s^3(c-4)+5cs^2-30s^2+20s-8cs-2c+24}{s(s+2)(cs^2-4s^2-c-s+6)}\ .$$
\end{theorem}

We use this formula to describe the cone $\mbox{Mov}(\mm_g)$ of
\emph{moving divisors}\footnote{Recall that an effective $\mathbb
Q$-Cartier divisor $D$ on a normal projective variety $X$ is said to
be moving, if the stable base locus $\bigcap_{n\geq 1}
\mbox{Bs}|\OO_X(nD)|$ has codimension at least $2$ in $X$.} inside
the cone $\mbox{Eff}(\mm_g)$ of effective divisors. The cone
$\mathrm{Mov}(\mm_g)$ parameterizes rational maps from $\mm_g$ in
the projective category while the cone $\mathrm{Nef}(\mm_g)$ of
numerically effective divisors, parameterizes regular maps from
$\mm_g$ (see \cite{HK} for details on this perspective). A
fundamental question is to estimate the following slope invariants
associated to $\mm_g$:
$$s(\mm_g):=\mathrm{inf}_{D\in \mathrm{Eff}(\mm_g)} \ s(D)$$  and
$$ s'(\mm_g):=\mathrm{inf}_{D\in \mathrm{Mov}(\mm_g)} \
s(D).$$ The formula  of the class of Brill-Noether divisors $\mm_{g,
d}^r$ when $\rho(g, r, d)=-1$  shows that
$\mathrm{lim}_{g\rightarrow \infty} s(\mm_g)\leq 6$ (cf.
\cite{EH2}). In \cite{F1} we provided an infinite sequence of genera
of the form  $g=a(2a+1)$ with $a\geq 2$ for which
$s(\mm_g)<6+12/(g+1)$, thus contradicting the Slope Conjecture
\cite{HMo}. There is no known example of a genus $g$ such that
$s(\mm_g)<6$.

Understanding the difference between $s(\mm_g)$ and $s'(\mm_g)$ is a
subtle question even for low $g$. There is a  strict inequality
$s(\mm_g)<s'(\mm_g)$\texttt{} whenever one can find an effective
divisor $D\in \mbox{Eff}(\mm_g)$ with $s(D)=s(\mm_g)$, such that
there exists a covering curve $R\subset D$ for which $R\cdot D<0$.
For $g<12$ the divisors minimizing the slope function have a strong
geometric characterization in terms of Brill-Noether theory. Thus
computing $s'(\mm_g)$ becomes a problem in understanding the
geometry of Brill-Noether and Gieseker-Petri divisors on $\mm_g$. To
illustrate this point we give two examples (see Section 5 for
details): It is known that $s(\mm_3)=9$ and the minimum slope is
realized by the locus of hyperelliptic curves $\mm_{3, 2}^1\equiv 9
\lambda-\delta_0-3\delta_1$. However $[\mm_{3, 2}^1]\notin
\mathrm{Mov}(\mm_3)$, because $\mm_{3, 2}^1$ is swept out by pencils
$R\subset \mm_3$ with $R\cdot \delta/R\cdot \lambda=28/3>s(\mm_{3,
2}^1)$. In fact, one has equality $s'(\mm_3)=28/3$ and the moving
divisor on $\mm_3$ attaining this bound corresponds to the pull-back
of an ample class under the rational map
$$\mm_3-->\mathcal{Q}_4:=|\OO_{\PP^2}(4)| \ //SL(3)$$ to the GIT
quotient of plane quartics which contracts $\mm_{3, 2}^1$ to a point
(see \cite{HL} for details on the role of this map in carrying out
the Minimal Model Program for $\mm_3$).

For $g=10$, it is known that $s(\mm_{10})=7$ and this bound is
attained by the divisor $\overline{\mathcal{K}}_{10}$ of curves
lying on $K3$ surfaces (cf. \cite{FP} Theorems 1.6 and 1.7 for
details):
$$\overline{\mathcal{K}}_{10}\equiv 7\lambda-\delta_0-5\delta_1-9 \delta_2-
12 \delta_3-14 \delta_4-15 \delta_5.$$ Furthermore,
$\overline{\mathcal{K}}_{10}$ is swept out by pencils $R\subset
\mm_{10}$ with $R\cdot \delta/R\cdot
\lambda=78/11>s(\overline{\mathcal{K}}_{10})$ (cf. \cite{FP}
Proposition 2.2). Therefore $[\overline{\mathcal{K}}_{10}]\notin
\mathrm{Mov}(\mm_{10})$ and $s'(\mm_{10})\geq 78/11$.

For $g=2s+1$ using the elementary observation that
$\phi^*(\mathrm{Ample}(\mm_{g'}))\subset \mathrm{Mov}(\mm_g)$,
Theorem \ref{slopepull} provides a uniform upper bound on slopes of
moving divisors on $\mm_g$:

\begin{corollary}\label{asy}
We set $g:=2s+1$ as and $g':=1+\frac{s}{s+1}{2s+2\choose s}$
as above. Then
$$ s(\phi^*(D))<6+\frac{16}{g-1} \mbox{  }  \mbox{ for every divisor } D\in \mathrm{Ample}(\mm_{g'}).$$
In particular one has the estimate $s'(\mm_{g})< 6+ 16/(g-1)$, for every odd integer $g\geq 3.$
\end{corollary}
Since we also know that $\mathrm{lim}_{g\rightarrow \infty}
s(\mm_g)\leq 6$, Corollary \ref{asy} indicates that (at least
asymptotically, for large $g$) we cannot distinguish between
effective and moving divisors on $\mm_g$ . We ask whether it is true
that $\mathrm{lim}_{g\rightarrow \infty}
s(\mm_g)=\mathrm{lim}_{g\rightarrow \infty} s'(\mm_g)$?

At the heart of the description in codimension $1$ of the map
$\phi:\mm_g-->\mm_{g'}$ lies the computation of the cohomology class
of the compactified Gieseker-Petri divisor
$\overline{\mathcal{GP}}_{g, d}^r\subset \mm_g$  in the case when
$\rho(g, r, d)=1$. Since this calculation is of independent interest
we discuss it in some detail. We denote by $\mathfrak{G}^r_d$ the
stack parameterizing pairs $[C, l]$ with $[C]\in \cM_g$ and $l=(L,
V)\in G^r_d(C)$ and denote by $\sigma:\mathfrak{G}^r_d\rightarrow
\cM_g$ the natural projection. In \cite{F1} we computed the class of
$\gp_{g, d}^r$ in the case $\rho(g, r, d)=0$, when $\gp_{g, d}^r$
can be realized as the push-forward of a determinantal divisor on
$\mathfrak{G}_d^r$ under the generically finite map $\sigma$. In
particular, we showed that if we write $g=rs+s$ and $d=rs+r$ where
$r\geq 1$ and  $s\geq 2$ (hence $\rho(g, r, d)=0)$, then we have the
following formula for the slope of $\gp_{g, d}^r$ (cf. \cite{F1},
Theorem 1.6):
$$s(\overline{\mathcal{GP}}_{g, d}^r)=6+ \frac{12}{g+1}+\frac{6(s+r+1)(rs+s-2)(rs+s-1)}{s(s+1)(r+1)(r+2)(rs+s+4)(rs+s+1)}.$$
The number $6+12/(g+1)$ is the slope of all Brill-Noether divisors
on $\mm_g$, that is $s(\gp_{g, d}^r)=6+12/(g+1)$ whenever $\rho(g,
r, d)=-1$ (cf. \cite{EH2}, or \cite{F1} Corollary 1.2 for a
different proof, making use of M. Green's Conjecture on syzygies of
canonical curves).

In the technically much-more intricate case $\rho(g, r, d)=1$, we
can realize $\cgp_{g, d}^r$ as the push-forward of a codimension $2$
determinantal subvariety of $\mathfrak{G}_d^r$ and most of Section 2
is devoted to extending this structure over a partial
compactification of $\mm_g$ corresponding to tree-like curves. If
$\sigma:\widetilde{\mathfrak G}^r_d \rightarrow \widetilde{\cM}_g$
denotes the stack of limit linear series $\mathfrak g^r_d$, we
construct two \emph{locally free} sheaves $\F$ and $\cN$ over
$\widetilde{\mathfrak G}^r_d$ such that $\mbox{rank}(\F)=r+1$,
$\mbox{rank}(\cN)=g-d+r=:s$ respectively, together with a vector bundle morphism
$$\mu: \F\otimes \cN\rightarrow \sigma^*\bigl(\mathbb E\otimes
\OO_{\mm_g}(\sum_{j=1}^{[g/2]} (2j-1)\cdot \delta_j)\bigr)$$ such
that $\gp_{g, d}^r$ is the push-forward of the first degeneration locus of
$\mu$:
\begin{theorem}\label{gisp2}
We fix integers $r, s\geq 1$ and we set $g:=rs+s+1, d:=rs+r+1$ so
that $\rho(g, r, d)=1$. Then the class of the compactified
Gieseker-Petri divisor $\gp_{g, d}^r$ in $\mm_g$ is given by the
formula:
$$\overline{\mathcal{GP}}^r_{g, d}\equiv \frac{C_{r+1} \ (s-1)r}{2(r+s+1)(s+r)(r+s+2)(rs+s-1)}\bigl(a \lambda-b_0 \delta_0-b_1
\delta_1 -\sum_{j=2}^{[g/2]} b_j \delta_j\bigr), $$ where
$$C_{r+1}:=\frac{(rs+s)!\ r! \ (r-1)! \cdots  2! \ 1!}{(s+r)!\ (s+r-1)! \cdots (s+1)! \ s!}$$
$$a=2s^3(s+1)r^5+s^2(2s^3+14s^2+33s+25)r^4+s(10s^4+ 59s^3+ 162s^2+179s+54)r^3+$$
$$+(18s^5+138s^4+387s^3+491s^2+244s+24)r^2+(14s^5+145s^4+464s^3+627s^2+378s+72)r+$$
$$4s^5+54s^4+208s^3+314s^2+212s+48  $$
$$b_0:= \frac{(r+2)(s+1)(s+r+1)(2rs+2s+1)(rs+s+2)(rs+s+6)}{6}  $$
$$b_1:=(r+1)s\Bigl(2s^2(s+1)r^4+s(2s^3+12s^2+23s+9s)r^3+(8s^4+39s^3+75s^2+46s+10)r^2+
$$
$$+(10s^4+59s^3+108s^2+89s+26)r+ 4s^4+30s^3+64s^2+58s+12\Bigr),$$
and $ b_j\geq b_1$ for $ j\geq 2$ are explicitly determined constants.
\end{theorem}
Even though the coefficients $a$ and $b_1$ look rather unwieldy, the
expression for the slope of $\overline{\mathcal{GP}}_{g, d}^r$ has a
simpler and much more suggestive expression which we record:
\begin{corollary}\label{slopegis}
For $\rho(g, r, d)=1$, the slope of the Gieseker-Petri divisor
$\gp_{g, d}^r$ has the following expression:
$$s(\overline{\mathcal{GP}}_{g, d}^r)=6+ \frac{12}{g+1}+\frac{24\ s(r+1)(r+s)(s+r+2)(rs+s-1)}
{(r+2)(s+1)(s+r+1)(2rs+2s+1)(rs+s+2)(rs+s+6)}.$$
\end{corollary}
Next we specialize to the case $r=1$, thus $g=2s+1$. Using the base point free pencil trick one can see
that the divisor $\cgp_{2s+1, s+2}^1$ splits into two irreducible components according to whether the pencil for which the
 Gieseker-Petri theorem fails has a base point or not. Precisely we have the following equality of codimension $1$ cycles
 $$\gp_{2s+1, s+2}^1=(2s-2)\cdot \mm_{2s+1, s+1}^{1}+\gp_{2s+1, s+2}^{1, 0},$$ where $\gp_{2s+1, s+2}^{1, 0}$ is the
 closure of the locus of curves $[C]\in \cM_g$ carrying a base point free pencil  $L\in W^1_{s+2}(C)$ such
 that $\mu_0(L)$ is not injective. Since we also have the well-known formula for the class of the Hurwitz divisor
 (cf. \cite{EH2}, Theorem 1)
$$\mm_{2s+1, s+1}^1\equiv \frac{(2s-2)!}{(s+1)! \ (s-1)!}\Bigl(6(s+2) \lambda-(s+1) \delta_0-6s \delta_1
-\cdots\Bigr),$$ we find the following expression for the slope of
$\gp_{2s+1, s+2}^{1, 0}$:
\begin{corollary}\label{slopepencil}
For $g=2s+1$, the slope of the divisor $\gp_{2s+1, s+2}^{1, 0}$ of
curves carrying a base point free pencil $L\in W^1_{s+2}(C)$ such
that $\mu_0(L)$ is not injective, is given by the formula
$$s( \gp_{2s+1, s+2}^{1, 0})=6+\frac{12}{g+1}+\frac{2s-1}{(s+1)(s+2)}.$$
\end{corollary}
We note that for $s=2$ and $g=5$, the divisor $\gp_{5, 4}^{1, 0}$ is equal to Teixidor's divisor of curves
 $[C]\in \cM_5$ having a vanishing theta-null, that is, a theta-characteristic
 $L^{\otimes 2}=K_C$ with $h^0(C, L)\geq 2$. In this case
 Corollary \ref{slopepencil} specializes to her formula
\cite{T} Theorem 3.1: $$\gp_{5, 4}^{1, 0}\equiv 4\cdot(33 \lambda-4
\delta_0-15 \delta_1-21  \delta_2)\in \mbox{Pic}(\mm_5).$$

To give another example we specialize to the case $r=1$, $s=3$ when
$g=7$. Using the base point free pencil trick, the divisor $\gp_{7,
5}^1$ can be identified with the closure of the locus of curves
$[C]\in \cM_7$ possessing a linear series $l\in G^2_7(C)$ such that
the plane model $C\stackrel{l}\rightarrow \PP^2$ has $8$ nodes, of
which $7$ lie on a conic. Its class is given by the formula:
$$\gp_{7, 5}^1\equiv 4\cdot (201  \lambda-26  \delta_0- 111
\delta_1-177  \delta_2- 198  \delta_3)\in \mathrm{Pic}(\mm_7).$$

In Section 5 we shall need a characterization of the $k$-gonal loci $\mm_{g, k}^1$ in terms of effective divisors of $\mm_g$ containing them. For instance, it is known that if $D\in \mathrm{Eff}(\mm_g)$ is a divisor such that $s(D)<8+4/g$, then $D$ contains the hyperelliptic locus $\mm_{g, 2}^1$ (see e.g. \cite{HMo}, Corollary 3.30). Similar bounds exist for the trigonal locus: if $s(D)<7+6/g$ then $D\supset \mm_{g, 3}^1$.  We have the following extension of this type of result:

\begin{theorem}\label{4gonal1}
1) Every effective divisor $D\in \mathrm{Eff}(\mm_g)$ having slope $$s(D)<\frac{1}{g} \Bigl[\frac{13g+16}{2}\Bigr]$$ contains the locus  $\mm_{g, 4}^1$ of $4$-gonal curves.

\noindent 2) Every effective divisor $D\in \mathrm{Eff}(\mm_g)$ having slope $$s(D)<\frac{1}{g}\Bigl(5g+9+2\Bigl[\frac{g+1}{2}\Bigr]\Bigr)$$ contains the locus $\mm_{g, 5}^1$ of $5$-gonal curves.
\end{theorem}
The proof uses an explicit unirational parametrization of $\mm_{g, k}^1$ that is available only when $k\leq 5$. It is natural to ask whether the subvariety $\mm_{g, k}^1\subset \mm_g$ is cut out by divisors $D\in \mathrm{Eff}(\mm_g)$ of slope less than the bound given in Theorem \ref{4gonal1}. Very little seems to be known about this question even in the hyperelliptic case.

We close by summarizing the structure of the paper. In Section 1 we
introduce a certain stack of pairs of complementary limit linear
series which we then use to prove Theorem \ref{gisp1} by induction on the
genus. The class of the compactified Gieseker-Petri divisor is computed in
Section 2. This calculation is  used in Section 3 to describe maps
between moduli spaces of curves. We then study
the geometry of $\phi$ in low genus (Section 4) with applications to
Prym varieties and we finish the paper by computing the invariant $s'(\mm_g)$
for $g\leq 11$ (Section 5).

\section{Divisorial components of the Gieseker-Petri locus}

Let us fix positive integers $g, r$ and $g$ such that $\rho(g, r,
d)\geq 0$ and set $s:=g-d+r\geq 2$, hence $g=rs+s+j$ and $d=rs+r+j$,
with $j\geq 0$. The case $j=0$ corresponds to the situation $\rho(g,
r, d)=0$ when we already know that $\cgp_{g, d}^r$ has a divisorial
component in $\cM_g$ whose class has been computed (see \cite{F1},
Theorem 1.6). We present an inductive method on $j$ which produces a
divisorial component of $\cgp_{g, d}^r \subset \cM_g$ provided one
knows that $\cgp_{g-1, d-1}^r$ has a divisorial component in
$\cM_{g-1}$. The method is based on degeneration to the boundary
divisor $\Delta_1\subset \mm_g$ and is somewhat similar to the one
used in \cite{F2} for the case $r=1$.

We briefly recall a few facts about (degeneration of) multiplication
maps on curves. If $L$ and $M$ are line bundles on a smooth curve
$C$, we denote by $$\mu_0(L, M): H^0(L)\otimes H^0(M)\rightarrow
H^0(L\otimes M)$$ the usual multiplication map and by $$\mu_1(L, M):
\mbox{Ker}\ \mu_0(L, M)\rightarrow H^0(K_C\otimes L\otimes M),\
\mbox{ }\mbox{ } \mu_1( \sum_{i} \sigma_i \otimes \tau_i):= \sum_i
 (d\sigma_i)\cdot \tau_i,$$ the first Gaussian map associated to
$L$ and $M$ (see \cite{W}). For any $\rho \in H^0(L)\otimes H^0(M)$
and a point $p\in C$, we write that $\mbox{ord}_p(\rho)\geq k$, if
$\rho$ lies in the span of elements of the form $\sigma\otimes
\tau$, where $\sigma \in H^0(L)$ and $\tau \in H^0(M)$ are such that
$\mbox{ord}_p(\sigma)+\mbox{ord}_p(\tau)\geq k$. When $i=0, 1$, the
condition $\mbox{ord}_p(\rho)\geq i+1$ for a generic point $p\in C$,
is clearly equivalent to $\rho\in \mbox{Ker}\ \mu_i(L, M)$.

If $X$ is a tree-like curve and $l$ is a limit $\mathfrak g^r_d$ on
$X$, for an irreducible component $Y\subset X$ we denote by
$l_Y=(L_Y, V_Y\subset H^0(L_Y))$ the $Y$-aspect of $l$. For $p\in Y$
we denote by $\{a^{l_Y}_i(p)\}_{i=0\ldots r}$ the \emph{vanishing
sequence} of $l$ at $p$ and by $$\rho(l_Y, p):=\rho(g(Y), r,
d)-\sum_{i=0}^r (a^{l_Y}_i(p)-i)$$ the \emph{adjusted Brill-Noether
number} with respect to the point $p$ (see  \cite{EH1} for a general
reference on limit linear series).

We shall repeatedly use the following elementary observation already
made in \cite{EH3} and used in \cite{F2}: Suppose
$\{\sigma_i\}\subset H^0(L)$ and $\{\tau_j\}\subset H^0(M)$ are
bases of global sections with the property that
$\mbox{ord}_p(\sigma_i)=a_i^L(p)$ and
$\mbox{ord}_p(\tau_j)=a_j^M(p)$ for all $i$ and $j$. Then if $\rho
\in \mbox{Ker}\ \mu_0(L,M))$, there must exist two pairs of integers
$(i_1,j_1)\neq (i_2,j_2)$ such that
$\mbox{ord}_p(\rho)=\mbox{ord}_p(\sigma_{i_1})+\mbox{ord}_p(\tau_{j_1})=\mbox{ord}_p(\sigma_{i_2})+
\mbox{ord}_p(\tau_{j_2})$.

A technical tool in the paper is the stack $\nu:
\widetilde{\mathcal{U}}_{g, d}^r\rightarrow \widetilde{\cM}_g$ of
pairs of complementary limit linear series defined over a partial
compactification of $\cM_g$ which will be defined below. Then
$\cgp_{g, d}^r$ is the push-forward under $\nu_{| \nu^{-1}(\cM_g)}$
of a degeneration locus inside $\widetilde{\mathcal{U}}_{g, d}^r$.
We denote by $\mathfrak {Pic}^d$ the degree $d$ Picard stack over
$\cM_g$, that is, the \'etale sheafification of the Picard functor,
and by $\mathbb E$ the Hodge bundle over $\mm_g$.  We consider
$\mathfrak G^r_d\subset \mathfrak{Pic}^d$ to be the stack
parameterizing pairs $[C, l]$ with $l=(L, V) \in G^r_d(C)$ and the
projection $\sigma:\mathfrak G^r_d\rightarrow \cM_g$.

We set $\Delta_0^0\subset \Delta_0\subset \mm_g$ to be the locus of
curves $[C/y\sim q]$, where $[C, q]\in \cM_{g-1, 1}$ is
Brill-Noether general and $y \in C$ is an arbitrary point, as well
as their degenerations $[C\cup_q E_{\infty}]$, where $E_{\infty}$ is
a rational nodal curve, that is, $j(E_{\infty})=\infty$. For $1\leq
i\leq [g/2]$, we denote by $\Delta_i^0\subset \Delta_i$ the open
subset consisting of unions $[C\cup_y D]$, where $[C]\in \cM_i$ and
$[D, y]\in \cM_{g-i, 1}$ are Brill-Noether general curves but the
point $y\in C$ is arbitrary. Then if we denote by
$$\widetilde{\cM}_g:=\cM_g\cup \bigl(\cup_{i=0}^{[g/2]}
\Delta_i^0\bigr),$$ one can extend the covering $\sigma:\mathfrak
G^r_d\rightarrow \cM_g$ to a proper map $\sigma:\widetilde{\mathfrak
G}^r_d\rightarrow \widetilde{\cM}_g$ from the stack
$\widetilde{\mathfrak G}^r_d$ of limit linear series $\mathfrak
g^r_d$.

We now introduce the stack $\nu:\widetilde{\cU}_{g, d}^r\rightarrow
\widetilde{\cM}_g$ of complementary linear series: For $[C]\in
\cM_g$, the fibre $\nu^{-1}[C]$ parameterizes pairs $(l, m)$ where
$l=(L, V)\in G^r_d(C)$ and $m=(K_C\otimes L^{\vee}, W)\in
G^{g-d+r-1}_{2g-2-d}(C)$. If $[C=C_1\cup_y C_2]\in
\widetilde{\cM}_g$, where $[C_1, y]\in \cM_{i, 1}$ and $[C_2, y]\in
\cM_{g-i, 1}$, the fibre $\nu^{-1}[C]$ consists of pairs of limit
linear series $(l, m)$, where $l=\{(L_{C_1}, V_{C_1}), (L_{C_2},
V_{C_2})\}$ is a limit $\mathfrak g^r_d$ on $C$ and
$$m=\{\bigl(K_{C_1}\otimes \OO_{C_1}(2(g-i)\cdot p)\otimes L_{{C_1}}^{-1},
W_{C_1}\bigr), \bigl(K_{C_2}\otimes \OO_{C_2}(2i\cdot p)\otimes
L_{C_2}^{-1}, W_{C_2}\bigr)\}$$  is a limit $\mathfrak
g_{2g-2-d}^{g-d+r-1}$ on $C$ which is complementary to $l$. There is
a morphism of stacks $$\epsilon:\widetilde{\cU}_{g, d}^r\rightarrow
\widetilde{\mathfrak G}_{g, d}^r$$ which forgets the limit $\mathfrak
g_{2g-2-d}^{g-d+r-1}$ on each curve. Clearly $\sigma \circ \epsilon
=\nu$.

\begin{definition}
For a smooth curve $C$ of genus $g$, a Gieseker-Petri
\emph{$\mathfrak{(gp)}^r_{d}$-relation} consists of a pair of linear
series $(L,V)\in G^r_d(C)$ and $(K_C\otimes L^{\vee}, W)\in
G_{2g-2-d}^{g-d+r-1}(C)$, together with an element $\rho \in
\PP\mbox{Ker}\{\mu_0(V, W):V\otimes W\rightarrow H^0(K_C)\}.$

If $C=C_1\cup_p C_2$ is a curve of compact type with $C_1$ and $C_2$
being smooth curves with $g(C_1)=i$ and $g(C_2)=g-i$ respectively, a
\emph{$\mathfrak{(gp)}_{d}^r$-relation} on $C$ is a collection
$(l,m,\rho_1, \rho_2)$, where $[C, l, m]\in \widetilde{\cU}_{g,
d}^r$, and elements
$$\rho_1\in \PP\mbox{Ker}\{V_{C_1}\otimes W_{C_1}\rightarrow
H^0\bigl(K_{C_1}(2(g-i)p)\bigr)\},$$
$$ \rho_{2}\in \PP
\mbox{Ker}\{V_{C_2}\otimes W_{C_2}\rightarrow H^0\bigl(
K_{C_2}(2ip)\bigr)\}$$ satisfying the compatibility relation
$\mbox{ord}_p(\rho_1)+\mbox{ord}_p(\rho_2)\geq 2g-2$.
\end{definition}

For every curve $C$ of compact type, the variety
$\mathcal{Q}_d^r(C)$ of $\mathfrak{(gp)}_d^r$-relations has an
obvious determinantal scheme structure. One can construct a moduli
stack of $\mathfrak{(gp)}^r_d$-relations which has a natural
determinantal structure over the moduli stack of curves of compact
type. In particular one has a lower bound on the dimension of each
irreducible component of this space and we shall use this feature in
order to smooth $\mathfrak{(gp)}^r_d$-relations constructed over
curves from the divisor $\Delta_1$ to nearby smooth curves from
$\cM_g$. The proof of the following theorem is very similar to the
proof of Theorem 4.3 in \cite{F2} which dealt with the case $r=1$.
We omit the details.

\begin{theorem}\label{gprelations}
We fix integers $g, r, d$ such that $\rho(g, r, d)\geq 0$ and a
curve $[C:=C_1\cup_y C_2]\in \mm_g$ of compact type. We denote by
$\pi:\mathcal{C}\rightarrow B$ the versal deformation space of
$C=\pi^{-1}(0)$, with $0\in B$. Then there exists a quasi-projective
variety $\nu: \mathcal{Q}_d^r\rightarrow B$, compatible with base
change, such that the fibre over each point $b\in B$ parameterizes
$\mathfrak{(gp)}^r_d$-relations over $C_b$. Moreover, each
irreducible component of $\mathcal{Q}_d^r$ has dimension at least
$\mathrm{dim}(B)-1=3g-4$.
\end{theorem}

The dimensional estimate on $\mathcal{Q}^r_d$ comes from its
construction as a determinantal variety over $B$. Just like in the
case of $\widetilde{\cU}_{g, d}^r$, we denote by $\epsilon:
\mathcal{Q}^r_d\rightarrow \widetilde{\mathfrak G}^r_d$ the
forgetful map such that $\sigma\circ \epsilon=\nu$. We use the
existence of $\mathcal{Q}_d^r$ to prove the following inductive
result:

\begin{theorem} Fix integers $g, r, d$ such that $\rho(g, r, d)\geq
2$ and let us assume that $\cgp_{g, d}^r$ has a divisorial component
$\mathcal{D}$ in $\cM_g$ such that if $[C]\in \mathcal{D}$ is a
general point, then the variety $\mathcal{Q}_d^r(C)$ has at least
one $0$-dimensional component corresponding to two complementary
base point free linear series $(l, m)\in G^r_d(C)\times
G_{2g-2-d}^{g-d+r-1}(C)$, such that $[C, l]\in \widetilde{\mathfrak
G}^r_d$ is a smooth point. Then $\cgp_{g+1, d+1}^r$ has a divisorial
component $\mathcal{D}'$ in $\cM_{g+1}$ such that a general point
$[C']\in \mathcal{D}'$ enjoys the same properties, namely that
$\mathcal{Q}_{d+1}^r(C')$ possesses a $0$-dimensional component
corresponding to a pair of base point free complementary linear
series $(l', m')\in G^r_{d+1}(C')\times G_{2g-1-d}^{g-d+r-1}(C')$
such that  $[C', l']\in \widetilde{\mathfrak G}^{r}_{d+1}$ is a
smooth point.
\end{theorem}
\begin{proof} We choose  a general
curve $[C]\in \mathcal{D} \subset \cgp_{g, d}^r$, a general point
$p\in C$ and we set $[C_0:=C\cup_p E]\in \mm_{g+1}$, where $E$ is an
elliptic curve. By assumption, there exist base point free linear
series $l_0=(L ,V)\in G^r_{d}(C)$ and $m_0=(K_C\otimes L^{\vee},
W)\in G_{2g-2-d}^{s-1}(C)$, together with an element $\rho\in \PP
\mbox{Ker}\bigl(\mu_0(V, W)\bigr)$ such that $\mbox{dim}_{(l_0, m_0,
\rho)} \mathcal{Q}_{d}^r(C)=0$. In particular, then $\mbox{Ker}\
\mu_0(V, W)$ is $1$-dimensional. Let $\pi:\mathcal{C} \rightarrow B$
be the versal deformation space of $C_0=\pi^{-1}(0)$ and
$\Delta\subset B$ the boundary divisor corresponding to singular
curves. We consider the scheme $\nu:\mathcal{Q}_{d+1}^r\rightarrow
B$ parameterizing $\mathfrak{(gp)}_{d+1}^r$-relations (cf. Theorem
\ref{gprelations}). Since $[C, l_0]\in \mathfrak{G}^r_d$ is a smooth
point and $l_0$ is base point free, Lemma 2.5 from \cite{AC} implies
that $$\mu_1(V, W):\mbox{Ker}\ \mu_0(V, W)\rightarrow
H^0(K_C^{\otimes 2})$$ is injective, in particular $\mu_1(V,
W)(\rho)\neq 0$. (Here $\sigma_0:\mathfrak{G}^r_d\rightarrow \cM_g$
denotes the stack of $\mathfrak g^r_d$'s over the moduli space of
curves of genus $g$).
 Thus we can assume that $\mbox{ord}_p(\rho)=1$ for a generic choice
of $p$.

We construct a $\mathfrak{(gp)}_{d+1}^r$-relation $z=(l,m, \rho_C,
\rho_E)\in \mathcal{Q}_{d+1}^r(C_0)$ as follows: the $C$-aspect of
the limit $\mathfrak g^r_{d+1}$ denoted by $l$ is obtained by adding
$p$ as a base point to $(L, V)$, that is $l_C= \bigl(L_{C}:=L\otimes
\OO_C(p), V_{C}:=V\subset H^0(L_{C})\bigr)$. The aspect $l_E$ is
constructed by adding $(d-r)\cdot p$ as a base locus to $|L^0_E|$,
where $L^0_E\in \mbox{Pic}^{r+1}(E)$ is such that $L_E^0\neq
\OO_E((r+1)\cdot p)$ and $(L^0_E)^{ \otimes 2}=\OO_E((2r+2)\cdot
p)$, and where $|V_E|=(d-r)\cdot p+|L_E^0|$. Since $p\in C$ is
general, we may assume that $p$ is not a ramification point of
$l_0$, which implies that $a^{l_C}(p)=(1, 2, \ldots, r+1)$. Clearly,
$$a^{l_E}(p)=(d-r, d-r+1, \cdots,  d),$$ hence $l=\{l_C, l_E\}$ is a
refined limit $\mathfrak g^r_{d+1}$ on $C_0$. The $C$-aspect of the
limit $\mathfrak g_{2g-2-d}^{s-1}$ we denote by $m$, is given by
$$m_C:=\bigl(K_C\otimes L^{\vee}\otimes \OO_C(p), W_C:=W\subset
H^0(K_C\otimes L^{\vee}\otimes \OO_C(p))\bigr).$$ The aspect $m_E$ is
constructed by adding $(g-r-1)\cdot p$ to the complete linear series
$|\OO_E((r+1+s)\cdot p)\otimes (L^0_E)^{\vee}|$. Since we may also
assume that $p$ is not a ramification point of $m_0$, we find that
$a^{m_C}(p)=(1, 2, \ldots, s)$ and $a^{m_E}(p)=(g-r-1, g-r, \ldots,
2g-2-d)$, that is, $m=\{m_C, m_E\}$ is a refined limit $\mathfrak
g_{2g-1-d}^{s-1}$ on $C_0$. Next we construct the elements $\rho_C$
and $\rho_E$. We choose
$$\rho_C=\rho \in \PP\mbox{Ker}\{\mu_0(V, W): V\otimes W\rightarrow H^0(K_C\otimes
\OO_C(2p))\},$$ that is, $\rho_C$ equals $\rho$ except that we add
$p$ as a simple base point to both linear series $l_C$ and $m_C$
whose sections get multiplied. Clearly $\mbox{ord}_p(\rho_C)=
\mbox{ord}_p(\rho)+2=3$. Then we construct an element $$\rho_E\in
\PP\mbox{Ker}\{V_E\otimes W_E\rightarrow H^0(\OO_E(2g\cdot p))\}$$
with the property that $\mbox{ord}_p(\rho_E)= 2g-3\
\bigl(=d-1+(2g-2-d)=d+(2g-3-d) \bigr)$. Such an element lies
necessarily in the kernel of the map
$$H^0\bigl(L_E^0 \otimes \OO_E(-(r-1)\cdot p) \bigr)\otimes
H^0\bigl(\OO_E((r+3)\cdot p)\otimes (L_E^0)^{\vee}\bigr)\rightarrow
H^0(\OO_E(4\cdot p)),$$ which by the base point free pencil trick is
isomorphic to the $1$-dimensional space $H^0\bigl(E,
\OO_E((2r+2)\cdot p)\otimes (L_E^0)^{\otimes (-2)}\bigr)$, that is,
$\rho_E$ is uniquely determined by the property that
$\mbox{ord}_p(\rho_E)\geq 2g-3$.

Since $\mbox{ord}_p(\rho_C)+\mbox{ord}_p(\rho_E)= 2g$, we find that
$z=(l, m, \rho_C, \rho_E)\in \mathcal{Q}_{d+1}^r$. Theorem
\ref{gprelations} guarantees that any component of
$\mathcal{Q}_{d+1}^r$ passing through $z$ has dimension at least
$3g-1$. To prove the existence of a component of
$\mathcal{Q}_{d+1}^r$ mapping rationally onto a divisor
$\mathcal{D}' \subset \cM_{g+1}$, it suffices to show that $z$ is an
isolated point in $\nu^{-1}([C_0])$. Suppose that $z'=(l', m',
\rho_C', \rho_E')\in \mathcal{Q}_{d+1}^r$ is another point lying in
the same component of $\nu^{-1}([C_0])$ as $z$. Since the scheme
$\mathcal{Q}^r_{d+1}$ is constructed as a disjoint union over the
possibilities of the vanishing sequences of the limit linear series
$\mathfrak g^r_{d+1}$ and $\mathfrak g^{s-1}_{2g-1-d}$, we may
assume that $a^{l'_C}(p)=a^{l_C}(p)=(1, 2, \ldots, r+1)$,
$a^{m'_C}(p)=a^{m_C}(p)=(1, 2, \ldots, s)$. Similarly for the
$E$-aspects, we assume that $a^{l'_E}(p)=a^{l_E}(p)$ and
$a^{m'_E}(p)=a^{m_E}(p)$. Then necessarily, $\mbox{ord}_p(\rho_C')=3
(=1+2=2+1)$, otherwise we would contradict the assumption $\mu_1(V,
W)(\rho)=0$. Moreover, $l_C=l_0$ and $m_C=m_0$ because of the
inductive assumption on $[C]$. Using the compatibility relation
between $\rho'_C$ and $\rho'_E$ we then get that
$\mbox{ord}_p(\rho'_E)\geq 2g-3$. The only way this can be satisfied
is when the underlying line bundle $L'_E$ of the linear series
$l'_E(-(d-r)\cdot p)$ satisfies the relation $(L_E')^{\otimes
2}=\OO_E((2r+2)\cdot p)$, which gives a finite number of choices for
$l_E'$ and then for $m_E'$. Once $l_E'$ is fixed, then as pointed
out before, $\rho'_E$ is uniquely determined by the condition
$\mbox{ord}_p(\rho'_E)\geq 2g-3$ (and in fact one must have
equality). This shows that $z\in \nu^{-1}([C_0])$ is an isolated
point, thus $z$ must smooth to $\mathfrak{(gp)}_{d+1}^r$ relations
on smooth curves filling-up a divisor $\mathcal{D}'$ in $\cM_{g+1}$.

We now prove that $[C_0, l]\in \widetilde{\mathfrak G}^r_{d+1}$ is a
smooth point (Recall that $\sigma: \widetilde{\mathfrak
G}^r_{d+1}\rightarrow B$ denotes the stack of limit $\mathfrak
g^r_{d+1}$'s on the fibres of $\pi$). This follows once we show that
$[C_0, l]$ is a smooth point of $\sigma^*(\Delta)$ and then observe
that $\widetilde{\mathfrak G}^r_{d+1}$ commutes with base change. By
explicit description, a neighbourhood of $[C_0, l] \in
\sigma^*(\Delta)$ is locally isomorphic to an \'etale neighbourhood
of $(\mathfrak G^r_d\times _{\cM_g} \cM_{g, 1})\times \cM_{1, 1}$
around the point $\bigl([C, l_0], [C, y], [E, y]\bigr)$ and we can
use our inductive assumption that $\mathfrak{G}^r_d$ is smooth at
the point $[C, l_0]$.

Finally, we prove that a generic point $[C']\in \mathcal{D}'$
corresponds to a pair of base point free linear series $(l', m')\in
G^r_{d+1}(C')\times G^{s-1}_{2g-1-d}(C')$. Suppose this is not the
case and assume that, say, $l'\in G^r_{d+1}(C')$ has a base point.
As $[C', l']\in \widetilde{\mathfrak G}^r_{d+1}$ specializes to
$[C_0, l_0]$, the base point of $l'$ specializes to a point $y\in
(C_{0})_{\mathrm{reg}}$ (If the base point specialized to the $p\in
C\cap E$, then necessarily $l$ would be a non-refined limit
$\mathfrak g^r_{d+1}$). If $y\in C-\{p\}$ then it follows that
$l_0=l_C(-p)\in G^r_d(C)$ has a base point at $y$, which is a
contradiction. If $y\in E-\{p\}$, then $L_E^0$ must have a base
point at $y$ which is manifestly false.
\end{proof}

\section{The class of the Gieseker-Petri divisors.}
In this section we determine the class of the Gieseker-Petri divisor
$\gp_{g, d}^r$. We start by setting some notation. We fix integers
$r, s\geq 1$ and set $g:=rs+s+1$ and $d:=rs+r+1$, hence $\rho(g, r,
d)=1$. We denote by $\cM_g^0$ the open substack of $\cM_g$
consisting of curves $[C]\in \cM_g$ such that
$W^{r+1}_d(C)=\emptyset$. Since $\rho(g, r+1, d)=-r-s-1$, it follows
that $\mbox{codim}(\cM_g-\cM_g^0, \cM_g)\geq 3$.
 In this section we denote by
 $\mathfrak G^r_d\subset \mathfrak{Pic}^d$ the stack
parameterizing pairs $[C, l]$ with $[C]\in \cM_g^0$ and $l\in
G^r_d(C)$ and $\widetilde{\cM}_g :=\cM_g^0 \cup(\cup_{i=0}^{[g/2]}
\Delta_i^0)$. We have a natural projection $\sigma:\mathfrak
G^r_d\rightarrow \cM_g^0$. Furthermore, we denote by $\pi:\cM_{g,
1}^0\rightarrow \cM_g^0$ the universal curve and by $f:\cM_{g, 1}^0
\times_{\cM_{g}^0} \mathfrak G^r_d\rightarrow \mathfrak G^r_d$ the
second projection. Note that the forgetful map $\epsilon: \cU_{g,
d}^r\rightarrow \mathfrak G^r_d$ is an isomorphism over $\cM_g^0$,
and we make the identification between $\cU_{g, d}^r$ and $\mathfrak
G^r_d$ (This identification obviously no longer holds over
$\widetilde{\cM}_g-\cM_g^0$).

\noindent From general Brill-Noether theory it follows that there exists a
unique component of $\mathfrak{G}^r_d$ which maps onto $\cM_g^0$.
Moreover, any irreducible component $\mathcal{Z}$ of $\mathfrak
G^r_d$ of dimension $> 3g-3+\rho(g, r, d)$ has the property that
$\mbox{codim}\bigl(\sigma(\mathcal{Z}), \cM_g^0\bigr)\geq 2$ (see
\cite{F1}, Corollary 2.5 for a similar statement when $\rho(g, r,
d)=0$, the proof remains essentially the same in the case $\rho(g,
r, d)=1$).

If $\L$ is a Poincar\'e bundle over $\cM_{g, 1}^0\times_{\cM_{g}^0}
\mathfrak G^r_d$ (one may have to make an \'etale base change
$\Sigma\rightarrow \mathfrak G^r_d$ to ensure the existence of $\L$,
see \cite{Est}), we set $\F:=f_*(\L)$ and $\cN:=R^1 f_*(\L)$. By
Grauert's theorem, both $\F$ and $\cN$ are vector bundles over
$\mathfrak G^r_d=\cU_{g, d}^r$ with $\mbox{rank}(\F)=r+1$ and
$\mbox{rank}(\cN)=s$ respectively, and there exists a bundle
morphism $\mu: \F\otimes \cN\rightarrow \sigma^*(\mathbb E)$, which
over each point $[C, L]\in \mathfrak G^r_d$ restricts to the Petri
map $\mu_0(L)$. If $\cU:=Z_{rs+s-1}(\mu)$ is the first degeneration
locus of $\mu$, then clearly $\cgp_{g, d}^r=\sigma_*(\cU)$. Each
irreducible component of $\cU$ has codimension at most $2$ inside
$\mathfrak G^r_d$. We shall prove that every such component mapping
onto a divisor in $\cM_g$ is in fact of codimension $2$ (see
Proposition \ref{divisor}), which will enable us to use Porteous'
formula to compute its class. While the construction of $\F$
and $\cN$ clearly depends on the choice of the Poincar\'e bundle
$\mathcal{L}$ (and of $\Sigma$), it is easy to check that the degeneracy class $Z_{rs+s-1}(\mu)\in A^2(\mathfrak G^r_d)$ is independent of
such choices.

Like in \cite{F1}, our technique for determining the class of the
divisor $\gp_{g, d}^r$ is to intersect $\cU$  with pull-backs of
test curves sitting in the boundary of $\mm_{g}$: We fix a general
pointed curve $[C, q]\in \cM_{g-1, 1}$ and a general elliptic curve
$[E, y]\in \cM_{1, 1}$. Then we define the families
$$C^0:=\{C/y\sim q: y\in C\}\subset \Delta_0\subset \mm_{g} \mbox{ and }
C^1:=\{C\cup_y E:y \in C\}\subset \Delta_1\subset \mm_{g}.
$$
These curves intersect the generators of $\mbox{Pic}(\mm_{g})$ as
follows:
$$C^0\cdot \lambda=0,\  C^0\cdot \delta_0=-2g+2,\ C^0\cdot \delta_1=1,
\mbox{  } C^0\cdot \delta_j=0 \mbox{ for }2\leq j\leq  [g/2]$$
and
$$C^1\cdot \lambda=0, \ C^1\cdot \delta_0=0, \ C^1\cdot \delta_1=-2g+4, \mbox{  } C^1\cdot
\delta_j=0 \mbox{ for }2\leq j\leq [g/2].$$ Next we fix a genus
$[g/2]\leq j\leq g-2$ and general curves $[C]\in \cM_j, [D, y]\in
\cM_{g-j, 1}$.  We define the $1$-parameter family $C^j:=\{C^j_y=
C\cup_y D\}_{y\in C}\subset \Delta_j\subset \mm_{g}$. We have the
formulas
$$C^j\cdot \lambda =0,\mbox{ } \ C^j\cdot \delta_j=-2j+2 \mbox{ and
} C^j\cdot \delta_i=0\ \mbox{ for } i\neq j.$$ To understand the
intersections $C^j\cdot \gp_{g, d}^r$ for $0\leq j\leq [g/2]$, we
shall extend the vector bundles $\F$ and $\cN$ over the partial
compactification $\widetilde{\cU}_{g, d}^r$ constructed in Section
1.

The following propositions describe the pull-back surfaces
$\sigma^*(C^j)$ inside $\widetilde{\mathfrak G}^r_d$:

\begin{proposition}\label{limitlin1}
We set $g:=rs+s+1$ and fix general curves $[C]\in \cM_{rs+s}$ and
$[E, y]\in \cM_{1, 1}$ and consider the associated test curve
$C^1\subset \Delta_1\subset \mm_g$. Then we have the following
equality of $2$-cycles in $\widetilde{\mathfrak{G}}^r_d$:
$$\sigma^*(C^1)=X+X_1\times X_2+\Gamma_0\times Z_0+n_1\cdot Z_1+n_2\cdot Z_2+n_3\cdot Z_3,$$
where
$$X:=\{(y, L)\in C\times W^r_d(C): h^0(C, L\otimes \OO_C(-2y))=r\}$$
$$X_1:=\{(y, L)\in C\times W^r_d(C): h^0( L\otimes \OO_C(-2\cdot y))=r,
\ h^0(L\otimes \OO_C(-(r+2)\cdot y))=1\}$$
$$X_2:=\{(y, l)\in G^r_{r+2}(E): a^l_1(y)\geq 2, \ a^l_{r}(y)\geq
r+2\}\cong \PP\Bigl(\frac{H^0(\OO_E((r+2)\cdot y))}{H^0(\OO_E(r\cdot
y))}\Bigr)$$
$$\Gamma_0:=\{(y, A\otimes \OO_C(y)): y\in C, A\in
W^r_{d-1}(C)\} \mbox{ }, \ Z_0=G^r_{r+1}(E)=\mathrm{Pic}^{r+1}(E)$$
$$Z_1:=\{l\in G^r_{r+3}(E): a_1^l(y)\geq 3, \ a^l_r(y)\geq
r+3\}\cong \PP\Bigl(\frac{H^0(\OO_E((r+3)\cdot y))}{H^0(\OO_E(r\cdot
y))}\Bigr)$$
$$Z_2:=\{l\in G^r_{r+2}(E):a_2^l(y)\geq 3, \ a_r^l(y)\geq r+2\}\cong
\PP\Bigl(\frac{H^0(\OO_E((r+2)\cdot y))}{H^0(\OO_E((r-1)\cdot
y))}\Bigr)$$
$$Z_3:=\{l\in G^r_{r+2}(E): a_1^l(y)\geq 2\}=\bigcup_{z\in E} \PP\Bigl(\frac
{H^0(\OO_E((r+1)\cdot y+z))}{H^0(\OO_E((r-1)\cdot y+z))}\Bigr),$$
where the constants $n_1, n_2, n_3$ are explicitly known positive
integers.
\end{proposition}
\begin{proof}
Every point in $\sigma^*(C^1)$ corresponds to a limit $\mathfrak
g^r_d$, say $l=\{l_C, l_E\}$, on some curve $[C^1_y:= C\cup_y E]\in C^1$. By
investigating the possible ways of distributing the Brill-Noether
numbers $\rho(l_C, y)$ and $\rho(l_E, y)$ in a way such that the inequality
$1=\rho(g, r, d) \geq \rho(l_C, y)+\rho(l_E, y)$ is satisfied, we
arrive at the six components in the statement (We always use the elementary inequality
$\rho(l_E, y)\geq 0$, hence $\rho(l_C, y)\leq 1$). We mention that
$X$ corresponds to the case when $\rho(l_C, y)=1, \rho(l_E, y)=0$,
the surfaces $X_1\times X_2$ and $\Gamma_0\times Z_0$ correspond to
the case $\rho(l_C, y)=0, \rho(l_E, y)=0$, while $Z_1, Z_2, Z_3$
appear in the cases when $\rho(l_C, y)=-1, \rho(l_E, y)=1$. The
constants $n_i$ for $1\leq i\leq 3$ have a clear enumerative
meaning: First, $n_1$ is the number of points $y\in C$ for which
there exists $L\in W^r_d(C)$ such that $a^L(y)=(0, 2, 3, \ldots, r,
r+3)$. Then $n_2$ is the number of points $y\in C$ for which there
exists $L\in W^r_d(C)$ such that $a^L(y)=(0, 2, 3, \ldots, r-1, r+1,
r+2)$. Finally, $n_3$ is the number of points $y\in C$ which appear
as  ramification points for one of the finitely many linear series
$A\in W^r_{d-1}(C)$.
\end{proof}

Next we describe $\sigma^*(C^0)$ and we start by fixing more
notation. We choose a general pointed curve $[C, q]\in \cM_{rs+s,
1}$ and denote by $Y$ the following surface:
$$Y:=\{(y, L)\in C\times W^r_{d}(C): h^0(C, L\otimes
\OO_C(-y-q))=r\}.$$ Let $\pi_1:Y \rightarrow C$ denote the first
projection. Inside $Y$ we consider two curves corresponding to
$\mathfrak g^r_{d}$'s with a base point at $q$:
$$\Gamma_1:=\{(y, A\otimes \OO_C(y)): y\in C, A\in W^r_{d-1}(C)\} \ \mbox{ } \mbox{ and }$$
 $$\Gamma_2:=\{(y, A\otimes \OO_C(q)): y\in C, A\in
W^r_{d-1}(C)\},$$ intersecting transversally in
$n_0:=\#(W^r_{d-1}(C))$ points. Note that $\rho(g, r-1, d)=0$ and
$W^r_{d-1}(C)$ is a reduced $0$-dimensional cycle. We denote by $Y'$
the blow-up of $Y$
 at these $n_0$
points and at the points $(q, B)\in Y$ where $B\in W^r_{d}(C)$ is a
linear series with the property that $h^0(C, B\otimes
\OO_C(-(r+2)\cdot q))\geq 1$. We denote by $E_A, E_B\subset Y'$ the
exceptional divisors corresponding to $(q, A\otimes \OO_C(q))$ and
$(q, B)$ respectively, by $\epsilon: Y'\rightarrow Y$ the projection
and by $\widetilde{\Gamma}_1, \widetilde{\Gamma}_2\subset Y'$ the
strict transforms of $\Gamma_1$ and $\Gamma_2$.

\begin{proposition}\label{limitlin0}
Fix a general curve $[C, q]\in \cM_{rs+s, 1}$ and consider the
associated test curve $C^0\subset \Delta_0\subset \mm_{rs+s+1}$.
Then we have the following equality of $2$-cycles in
$\widetilde{\mathfrak{G}}_{r}^d$:
$$\sigma^*(C^0)=Y'/\widetilde{\Gamma}_1\cong \widetilde{\Gamma}_2,$$
that is, $\sigma^*(C^0)$ can be naturally identified with the
surface obtained from $Y'$ by identifying the disjoint curves
$\widetilde{\Gamma}_1$ and $\widetilde{\Gamma}_2$ over each pair
$(y, A)\in C\times W^r_{d-1}(C)$.
\end{proposition}

\begin{proof}
We fix a point $y\in C-\{q\}$, denote by $[C_y^0:=C/ y\sim q]\in
\mm_{g}$, $\nu:C\rightarrow C^y_0$ the normalization map, and we
investigate the variety $\overline{W}^r_{d}(C_y^0)\subset
\overline{\mbox{Pic}}^{d}(C_y^0)$ of torsion-free sheaves $L$ on
$C_y^0$ with $\mbox{deg}(L)=d$ and $h^0(C_y^0, L)\geq r+1$. If $L\in
W^r_{d}(C_y^0)$, that is, $L$ is locally free, then $L$ is
determined by $\nu^*(L)\in W^r_{d}(C)$ which has the property that
$h^0(C, \nu^*L\otimes \OO_C(-y-q))=r$. However, the line bundles of
type $A\otimes \OO_C(y)$ or $A\otimes \OO_C(q)$ with $A\in
W^r_{d-1}(C)$, do not appear in this association even though they
have this property. They correspond to the situation when $L\in
\overline{W}_{d}^r(C_0^y)$ is not locally free, in which case
necessarily $L=\nu_*(A)$ for some $A\in W^r_{d-1}(C)$. Thus $Y\cap
\pi_1^{-1}(y)$ is the partial normalization of
$\overline{W}_{d}^r(C^0_y)$ at the $n_0$ points of the form
$\nu_*(A)$ with $A\in W^r_{d-1}(C)$. A special analysis is required
when $y=q$, that is, when $C_y^0$ degenerates to $C\cup _q
E_{\infty}$, where $E_{\infty}$ is a rational nodal cubic. If
$\{l_C, l_{E_{\infty}}\}\in \sigma^{-1}([C\cup_{q} E_{\infty}])$,
then an analysis along the lines of Theorem \ref{limitlin1} shows
that $\rho(l_C, q)\geq 0$ and $\rho(l_{E_{\infty}}, q)\leq 1$. Then
either $l_C$ has a base point at $q$ and then the underlying line
bundle of $l_C$ is of type $A\otimes \OO_C(q)$ while
$l_{E_{\infty}}(-(d-r-1)\cdot q)\in
\overline{W}_{r+1}^r(E_{\infty})$, or else, $a^{l_C}(q)=(0, 2, 3,
\ldots, r, r+2)$ and then $l_{E_{\infty}}(-(d-r-2)\cdot q)\in
\PP\bigl(H^0(\OO_{E_{\infty}}((r+2)\cdot
q))/H^0(\OO_{E_{\infty}}(r\cdot q))\bigr)\cong E_B$, where $B\in
W^r_{d}(C)$ is the underlying line bundle of $l_C$.
\end{proof}

We now show that every irreducible component of $\cU$ has the
expected dimension:
\begin{proposition}\label{divisor}
Every irreducible component $\mathcal{X}$ of $\cU$ having the
property that $\sigma(\mathcal{X})$ is a divisor in $\cM_g$ has
$\mathrm{codim}(\mathcal{X}, \mathfrak G^r_d)=2$.
\end{proposition}
\begin{proof} Suppose that $\mathcal{X}$ is an irreducible component
of $\cU$ satisfying (1) $\mbox{codim}(\mathcal{X}, \mathfrak
G^r_d)\leq 1$ and (2) $\mbox{codim}(\sigma(\mathcal{X}), \cM_g)=1$.
We write $D:=\overline{\sigma(\mathcal{X})}\subset \mm_g$ for the
closure of this divisor in $\mm_g$, and we express its class as
$D\equiv a\lambda-b_0\delta_0-b_1\delta_1-\cdots
-b_{[g/2]}\delta_{[g/2]} \in \mbox{Pic}(\mm_g)$. To reach a
contradiction, it suffices to show that $a=0$.

Keeping the notation from Propositions \ref{limitlin1} and
\ref{limitlin0}, we are going to show that $C^0\cap D=C^1\cap
D=\emptyset$ which implies that $b_0=b_1=0$. Then we shall show that
if $R\subset \mm_g$ denotes the pencil obtained by attaching to a
general pointed curve $[C, q]\in \cM_{rs+s, 1}$ at the fixed point
$q$, a pencil of plane cubics (i.e. an elliptic pencil of degree
$12$), then $R\cap D=\emptyset$. This implies the relation
$a-12b_0+b_1=0$ which of course yields that $a=0$.

We assume by contradiction that $C^1\cap D\neq \emptyset$. Then
there exists a point $y\in C$ and a limit $\mathfrak g^r_d$ on
$C^1_y:=C\cup_y E$, say $l=\{l_C, l_E\}$, such that if $L_C\in
W^r_d(C)$ denotes the underlying line bundle of $l_C$, then the
multiplication map $$\mu_0(L_C, y): H^0(L_C)\otimes H^0(K_C\otimes
L_C^{\vee}\otimes \OO_C(2y))\rightarrow H^0(K_C\otimes \OO_C(2y))$$
is not injective. We claim that this can happen only when $\rho(l_C,
y)=1$ and $\rho(l_E, y)=0$, that is, when $[C^1_y, l]\in X$ (we are
still using the notation from Proposition \ref{limitlin1}). Indeed,
assuming that $\rho(l_C, y)\leq 0$, there are two cases to consider.
Either $L_C$ has a base point at $y$ and then we can write
$L_C=A\otimes \OO_C(y)$ for $A\in W^r_{d-1}(C)$ and then we find that
$\mu_0(A)$ is not injective which contradicts the assumption that
$[C]\in \cM_{rs+s}$ is Petri general. Or $y\notin \mbox{Bs}|L_C|$
and then $a^{L_C}(y)\geq (0, 2, 3, \ldots, r, r+2)$. A degeneration
argument along the lines of \cite{F1} Proposition 3.2 shows that
$[C]$ can be chosen general enough such that every $L_C$ with this
property has $\mu_0(L_C, y)$ injective. Thus we may assume that
$\rho(l_C, y)=1$ and then $\mu_0(L_C, y)$ is not injective for
\emph{every} point $(y, L_C)$ belonging to an irreducible component
of the fibre $\pi_1^{-1}(y)\subset X$.

On the other hand, whenever one has an irreducible projective
variety $A\subset G^r_d(C)$ with $\mbox{dim}(A)\geq 1$ and a
Schubert index $\overline{\alpha}:=(0\leq \alpha_0\leq \ldots \leq
\alpha_r\leq d-r)$ such that $\alpha^l(y)\geq \overline{\alpha}$ for
all $l\in A$, there exists a Schubert index of the same type
$\overline{\beta} >\overline{\alpha}$, such that
$\alpha^{l_0}(y)\geq \overline{\beta}$ for a certain $l_0\in A$. In
our case, this implies that $\mu_0(L_C, y)$ is not injective for a
linear series $L_C\in W^r_d(C)$ such that either $a^{L_C}(y)\geq (0,
2, \ldots, r, r+2)$ (and this case has been dealt with before), or
$a^{L_C}(y)\geq (1, 2, \ldots, r+1)$. Then $L_C=A\otimes \OO_C(y)$
for $A\in W^r_{d-1}(C)$ and $\mu_0(A)$ is not injective. This
violates the assumption that $[C]\in \cM_{rs+s}$ is Petri general.
To prove  that $C^0\cap D=\emptyset $ we use the same principle in
the context of the explicit description of $\sigma^*(C^0)$ provided
by Proposition \ref{limitlin0}. Finally, to show that $R\cap
D=\emptyset$ it suffices to show that if $[C, q]\in \cM_{rs+s, 1}$
is sufficiently general, then $\mu_0(L_C, q)$ is injective for every
$(q, L_C)\in \pi_1^{-1}(q)$. This is the statement of Theorem
\ref{elltail}.
\end{proof}
We extend  $\F$ and $\cN$ as vector bundles over the stack
$\widetilde{\cU}_{g, d}^r$ of pairs of limit linear series. Note
that every irreducible component of $\widetilde{\cU}_{g, d}^r$ which
meets one of the test surfaces $\nu^*(C^j)$ has dimension $3g-2$.
This follows from an explicit description of $\nu^*(C^j)$ similar to
the one for $j=0, 1$ given in Propositions \ref{limitlin1} and
\ref{limitlin0}. Such a description, although straightforward, is
combinatorially involved (see \cite{F1} Proposition 2.4, for the
answer in the case $\rho(g, r, d)=0$). Since we are not going to
make direct use of it in this paper, we skip such details. Recall
that we denote by $\epsilon: \widetilde{\cU}_{g, d}^r\rightarrow
\widetilde{\mathfrak G}^r_d$ the forgetful map and $\nu=\sigma\circ
\epsilon$.
\begin{proposition}
There exist two vector bundles $\F$ and $\cN$  over
$\widetilde{\cU}^r_{g, d}$ with $\mathrm{rank}(\F)=r+1$ and
$\mathrm{rank}(\cN)=s$, together with a vector bundle morphism
$\mu:\F\otimes \cN\rightarrow \nu^*\bigl(\mathbb E\otimes
\sum_{j=1}^{[g/2]} (2j-1)\cdot \delta_j\bigr)$, such that the
following statements hold:
\begin{itemize}
\item For a point $[C, L]\in \mathfrak{G}^r_{d}=\cU_{g, d}^r$
 we have that $\F(C, L)=H^0(C, L)$,\mbox{ } $\cN(C, L)=H^0(C, K_C\otimes
L^{\vee})$ and $$\mu_0(C, L): H^0(C, L)\otimes H^0(C, K_C\otimes
L^{\vee})\rightarrow H^0(K_C)$$ is the Petri map.

\item For $t=\bigl[C\cup_y D, (l_C, l_D), (m_C, m_D)\bigr]\in \sigma^{-1}(\Delta_j^0)$,
with $[g/2]\leq j\leq g-1,\  [C, y]\in \cM_{j, 1}$,  $[D, y]\in
\cM_{g-j, 1}$ and $$l_C=(L_C, V_C)\in G^r_d(C),$$
$$m_C=(K_C\otimes L_C^{\vee}\otimes \OO_C(2(g-j)\cdot y), W_C)\in
G^{s-1}_{2g-2-d}(C),$$ we have that  $\F(t)=V_C$, $\cN(t)=W_C$ and
$$\mu(t)=\mu_0(V_C, W_C): V_C\otimes W_C \rightarrow
H^0\bigl(K_C\otimes \OO_C(2(g-j)\cdot y)\bigr).$$

\item Fix $t=[C_y^0:=C/y\sim q, L]\in \sigma^{-1}(\Delta_0^0)$, with $q,
y\in C$ and $L\in \overline{W}^r_{d}(C_y^0)$  such that $h^0(C,
\nu^*L\otimes \OO_C(-y-q))=r$. Here $\nu:C\rightarrow C_y^0$ is the
normalization map.

\noindent When $L$ is locally free, $$\F(t)=H^0(C, \nu^*L),$$
$$\cN(t)=H^0(C, K_C\otimes \nu^*L^{\vee}\otimes \OO_C(y+q))$$  and
$\phi(t)$ is the multiplication map
$$H^0(\nu^*L)\otimes H^0\bigl(K_C\otimes \nu^*L^{\vee}\otimes
\OO_C(y+q)\bigr)\rightarrow H^0(K_C\otimes \OO_C(y+q))=H^0(C_y^0,
\omega_{C_y^0}).$$
 In the case when $L$ is not locally free, that is, $L\in
\overline{W}_{d}^r(C_0^y)-W_{d}^r(C_0^y)$, then $L=\nu_*(A)$, where
$A\in W^r_{d-1}(C)$, and
$$\F(t)=H^0(A)=H^0(\nu_*A)$$ and  
$$\cN(t)=H^0(K_C\otimes A^{\vee}\otimes
\OO_C(y+q))=H^0(\omega_{C^0_y}\otimes \nu_*A^{\vee}).$$
\end{itemize}
\end{proposition}

Briefly stated, over each curve of compact type, the vector bundle
$\F$ (resp. $\cN$) retains the sections of the limit $\mathfrak
g^r_d$ (resp. $\mathfrak g^{s-1}_{2g-2-d}$) coming from the
component having the largest genus. The Gieseker-Petri theorem
ensures that the vector bundle morphism $\mu:\F\otimes
\cN\rightarrow \nu^*\bigl(\mathbb E\otimes \sum_{j=1}^{[g/2]}
(2j-1)\cdot \delta_j\bigr)$ is generically non-degenerate. Moreover,
$\nu_{|\nu^{-1}(\Delta_0^0)}$ and $\nu_{| \nu^{-1}(\Delta_1^0)}$ are
also generically-nondegenerate along each irreducible component (see
Theorem \ref{elltail}), hence one can write that
$$\nu_* \ c_1\bigl(\nu^*(\mathbb E\otimes \sum_{j=1}^{[g/2]}
(2j-1)\cdot \delta_j)-\F\otimes \cN\bigr)=[\gp_{g,
d}^r]+\sum_{j=2}^{[g/2]} e_j\cdot \delta_j,$$ where $e_j\geq 0$. We
can compute explicitly the left-hand-side of this formula and show
that the smallest boundary coefficient of $\nu_*
c_1\bigl(\nu^*(\mathbb E\otimes \sum_{j=1}^{[g/2]} (2j-1)\cdot
\delta_j)-\F\otimes \cN\bigr)$ is that corresponding to $\delta_0$.
Thus $$s([\gp_{g, d}^r])=s\bigl(\nu_* \ c_1(\nu^*(\mathbb E\otimes
\sum_{j=1}^{[g/2]} (2j-1)\cdot \delta_j)-\F\otimes \cN)\bigr).$$

Throughout the paper we use a few facts about intersection theory on
Jacobians which we briefly recall (see \cite{ACGH} for a general
reference). We fix integers $r, s\geq 1$ and set $g:=rs+s$ and
$d:=rs+r+1$.
 If $[C]\in \cM_{g}$ is a Brill-Noether general curve, we denote
by $\P$ a Poincar\'e bundle on $C\times \mbox{Pic}^d(C)$ and by
$\pi_1:C\times \mbox{Pic}^d(C)\rightarrow C$ and $\pi_2:C\times
\mbox{Pic}^d(C)\rightarrow \mbox{Pic}^d(C)$ the projections. We
define the cohomology class $\eta=\pi_1^*([point])\in H^2(C\times
\mbox{Pic}^d(C))$, and if $\delta_1,\ldots, \delta_{2g}\in H^1(C,
\mathbb Z)\cong H^1(\mbox{Pic}^d(C), \mathbb Z)$ is a symplectic
basis, then we set
$$\gamma:=-\sum_{\alpha=1}^g
\Bigl(\pi_1^*(\delta_{\alpha})\pi_2^*(\delta_{g+\alpha})-\pi_1^*(\delta_{g+\alpha})\pi_2^*(\delta_
{\alpha})\Bigr).$$ We have the formula $c_1(\P)=d\cdot \eta+\gamma,$
corresponding to the Hodge decomposition of $c_1(\P)$. We also
record that $\gamma^3=\gamma \eta=0$, $\eta^2=0$ and
$\gamma^2=-2\eta \pi_2^*(\theta)$. Since $W^{r+1}_d(C)=\emptyset$,
it follows that $W^r_d(C)$ is smooth of dimension $\rho(g, r, d)=r+1
$. Over $W^r_d(C)$ there is a  tautological rank $r+1$ vector bundle
$\mathcal{M}:=(\pi_2)_{*}(\mathcal{P}_{| C\times W^r_d(C)})$. The
Chern numbers of $\mathcal{M}$ can be computed using the Harris-Tu
formula (cf. \cite{HT}) as follows: We write $$\sum_{i=0}^{r+1}
c_i(\mathcal{M}^{\vee})=(1+x_1)\cdots (1+x_{r+1})$$ and then for
every class $\zeta \in H^*(\mbox{Pic}^d(C), \mathbb Z)$ one has the
following formula:
$$x_1^{i_1}\cdots x_{r+1}^{i_{r+1}}\
\zeta=\mbox{det}\Bigl(\frac{\theta^{g+r-d+i_j-j+l}}{(g+r-d+i_j-j+l)!}\Bigr)_{1\leq
j, l\leq r+1}\ \zeta.$$ If we use the expression of the Vandermonde
determinant, we get the identity
$$\mbox{det}\Bigl(\frac{1}{(a_j+l-1)!}\Bigr)_{1\leq j, l\leq
r+1}=\frac{\Pi_{ j>l}\ (a_l-a_j)}{\Pi_{j=1}^{r+1}\ (a_j+r)!},$$
which quickly leads to the following formula in $H^{2r+2}(W^r_d(C),
\mathbb Z)$:
\begin{equation}\label{master}
x_1^{i_1}\cdots x_{r+1}^{i_{r+1}}\cdot
\theta^{r+1-i_1-\cdots-i_{r+1}}=\frac{\Pi_{j>l}(i_l-i_j+j-l)}
{\Pi_{j=1}^{r+1} (s+r+i_j-j)!}\theta^g.
\end{equation}
By repeatedly applying (\ref{master}), we get all intersection
numbers on $W^r_d(C)$ we shall need. We define the integer
$$n_0=C_{r+1}:=\frac{(rs+s)!\ r! \ (r-1)! \cdots  2! \ 1!}{(s+r)!\ (s+r-1)! \cdots (s+1)! \ s!}\ =\#(W^r_{d-1}(C))$$ and we have the following formulas:
\begin{proposition}\label{intersection}
Let $C$ be a general curve of genus $rs+s$ and we set $d:=rs+r+1$.
We denote by $c_i:=c_i(\cM^{\vee})\in H^{2i}(W^r_d(C), \mathbb Z)$
the Chern classes of the dual of the tautological bundle on
$W^r_d(C)$. Then one has the following identities in $H^*(W^r_d(C),
\mathbb Z)$:
$$c_{r+1}=x_1 x_2\ldots x_{r+1}=C_{r+1},$$ 
$$c_r\cdot c_1=x_1 x_2\ldots x_{r+1}+x_1^2 x_2\ldots x_r$$
$$c_{r-1}\cdot c_2=x_1x_2\ldots x_{r+1}+x_1^2x_2\ldots x_r+x_1^2 x_2^2x_3\ldots x_{r-1}$$
$$c_{r-1}\cdot c_1^2=x_1x_2 \ldots x_{r+1}+2x_1^2 x_2 \ldots x_r+x_1^2 x_2^2 x_3\ldots x_{r-1}+
x_1^3 x_2 x_3 \ldots x_{r-1} $$
$$c_r\cdot \theta=x_1x_2\ldots x_r\cdot \theta=(r+1)s \ C_{r+1},$$
$$c_{r-1}\cdot c_1 \cdot \theta=x_1x_2 \ldots x_r\cdot \theta+x_1^2 x_2 +\ldots x_{r-1}\cdot \theta $$
$$c_{r-2}\cdot c_2\cdot \theta=x_1 x_2 \ldots x_{r}\cdot \theta + x_1^2 x_2 \ldots x_{r-1}\cdot \theta+x_1^2 x_2^2 x_3 \ldots x_{r-2}\cdot \theta $$
$$c_{r-2}\cdot c_1^2\cdot \theta=x_1 x_2 \ldots x_r\cdot \theta+2x_1^2x_2\ldots x_{r-1}\cdot \theta +
x_1^2 x_2^2 x_3 \ldots x_{r-2}\cdot \theta+x_1^3 x_2x_3 \ldots x_{r-2}\cdot \theta $$
$$c_{r-1}\cdot \theta^2=x_1 x_2 \ldots x_{r-1}\cdot \theta^2,\  \mbox{ }\  \mbox{ }c_{r-2}\cdot c_1\cdot \theta^2=x_1 x_2 \ldots x_{r-1} \cdot \theta^2 +x_1^2 x_2 \ldots x_{r-2}\cdot \theta^2 $$
\end{proposition}

Next we record the values of the monomials in the $x_i$'s and
$\theta$ that appeared in Proposition \ref{intersection}. The proof
amounts to a systematic application of formula (\ref{master}):
\begin{proposition}\label{intersection2}
We set $d:=rs+r+1$ and  write $c_t(\cM^{\vee})=(1+x_1)\cdots
(1+x_{r+1})$ as above.  Then one has the following identities in
$H^{2r+2}(W_d^r(C), \mathbb Z)$:
$$ x_1 x_2\ldots x_{r+1}=C_{r+1},$$
$$x_1^2 x_2^2x_3\ldots x_{r-1}=
\frac{s(s+1)(r+1)^2(r-2)(r+2)}{4(s+r)(s+r+1)}C_{r+1},$$
$$ x_1^2 x_2\ldots x_r=\frac{r(r+2)s}{s+r+1}C_{r+1}, $$
$$
 x_1^3 x_2 x_3 \ldots x_{r-1}=\frac{r(r-1)(r+2)(r+3)s(s+1)}{4(s+r+1)(s+r+2)}C_{r+1}$$
$$ x_1x_2\ldots x_r\cdot \theta=(r+1)sC_{r+1},$$
 $$\  x_1^2 x_2 \ldots x_{r-1}\cdot \theta=\frac{(s+1)(r-1)(r+2)}{2(s+r+1)}x_1x_2\ldots x_r\cdot \theta, $$
$$ x_1^2 x_2^2 x_3 \ldots x_{r-2}\cdot \theta=\frac{(r-3)(r+1)(r+2)r(s+1)(s+2)}{12(s+r+1)(s+r)}x_1x_2\ldots x_r\cdot \theta$$
$$ \ x_1^3 x_2x_3 \ldots x_{r-2}\cdot \theta=
\frac{(r+2)(r+3)(r-1)(r-2)(s+1)(s+2)}{12(s+r+1)(s+r+2)}x_1x_2\ldots
x_r\cdot \theta$$
$$ x_1 x_2 \ldots x_{r-1}\cdot \theta^2=\frac{r(r+1)s(s+1)}{s}C_{r+1} $$
$$ x_1^2 x_2 \ldots x_{r-2}\cdot \theta^2=\frac{(r+2)(r-2)(s+2)}{3(s+r+1)}x_1 x_2 \ldots x_{r-1}\cdot \theta^2$$
$$ x_1 x_2 \ldots x_{r-2}\cdot \theta^3=\frac{(r+1)r(r-1)(s+2)(s+1)s}{6}C_{r+1}.$$
\end{proposition}

\begin{proposition}\label{xy}
Let $[C, q]\in \cM_{rs+s, 1}$ be a  general pointed curve. If $\cM$
denotes the tautological vector bundle over $W^r_{d}(C)$ and
$c_i:=c_i(\cM^{\vee})$, then one has the following relations:
\begin{enumerate}
\item
$[X]=\pi_2^*(c_r)-6\pi_2^*(c_{r-2})\eta \theta+
\bigl((4rs+2r+2s)\eta+2\gamma\bigr) \pi_2^*(c_{r-1}) \in
H^{2r}(C\times W^r_{d}(C))$.
\item
$[Y]=\pi_2^*(c_r)-2\pi_2^*(c_{r-2})\eta \theta+
\bigl((rs+r)\eta+\gamma) \pi_2^*(c_{r-1}) \in H^{2r}(C\times
W^r_{d}(C))$.
\end{enumerate}
\end{proposition}
\begin{proof}
We realize the surface $X$ as the degeneracy locus of a vector
bundle map over $C\times W^r_{d}(C)$. For each pair $(y, L)\in
C\times W^r_{d}(C)$ there is a natural map
$$H^0(C, L\otimes \OO_{2y})^{\vee}\rightarrow H^0(C, L)^{\vee}$$
which globalizes to a vector bundle morphism $\zeta:
J_1(\mathcal{P})^{\vee} \rightarrow \pi_2^*(\cM)^{\vee}$ over
$C\times W^r_{d}(C)$ (Recall that  $W^r_{d}(C)$ is a smooth
$(r+1)$-fold). Then we have the identification $X=Z_1(\zeta)$ and
the Thom-Porteous formula gives that $[X]=c_r\bigl(\pi^*_2(\cM)^{\vee}-
J_1(\mathcal{P}^{\vee})\bigr).$ From the usual exact sequence over
$C\times \mbox{Pic}^{d}(C)$
$$
0\longrightarrow \pi_1^*(K_C)\otimes \mathcal{P} \longrightarrow
J_1(\mathcal{P}) \longrightarrow \mathcal{P} \longrightarrow 0, $$
we can compute the total Chern class of the jet bundle
$$c_t(J_1(\mathcal{P})^{\vee})^{-1}=\Bigl(\sum_{j\geq
0}(d\eta+\gamma)^j\Bigr)\cdot \Bigl(\sum_{j\geq
0}((2g(C)-2+d)\eta+\gamma)^j\Bigr)=1-6\eta \theta +(2d+2g(C)-2)\eta
+2\gamma,
$$ which quickly leads to the formula for $[X]$. To compute $[Y]$ we
proceed in a similar way. We denote by $p_1, p_2:C\times C\times
\mbox{Pic}^{d}(C)\rightarrow C\times \mbox{Pic}^{d}(C)$ the two
projections, by $\Delta\subset C\times C\times \mbox{Pic}^{d}(C)$
the diagonal  and we set $\Gamma_q:=\{q\}\times \mbox{Pic}^{d}(C)$.
We introduce the rank $2$ vector bundle
$\cB:=(p_1)_*\bigl(p_2^*(\mathcal{P})\otimes
\OO_{\Delta+p_2^*(\Gamma_q)}\bigr)$ defined over $C\times
W^r_{d}(C)$ and we note that there is a bundle morphism $\chi:
\cB^{\vee}\rightarrow (\pi_2)^*(\cM)^{\vee}$ such that
$Y=Z_1(\chi)$. Since we also have that
$$c_t(\cB^{\vee})^{-1}=\bigl(1+(d\eta+\gamma)+(d\eta+\gamma)^2+\cdots\bigr)(1-\eta),$$
we immediately obtained the desired expression for $[Y]$.
\end{proof}

\begin{remark}
For future reference we also record the following formulas:
\begin{equation}\label{c3d1}
c_{r+1}(\pi_2^*(\cM)^{\vee}-J_1(\P)^{\vee})=\pi_2^*(c_{r+1})-6
 \pi_2^*(c_{r-1})\eta \theta + \bigl((4rs+2r+2s)\eta+2\gamma \bigr)\pi_2^*(c_r)
\end{equation}
\begin{equation}\label{c4d1}
c_{r+2}(\pi_2^*(\cM)^{\vee}-J_1(\P)^{\vee})=\pi_2^*(c_{r+1})\bigl((4rs+2r+2s)
\eta+ 2\gamma\bigr)-6 \pi_2^*(c_r) \eta \theta.
\end{equation}
\begin{equation}\label{c3d2}
c_{r+1}(\pi_2^*(\cM)^{\vee}-\cB^{\vee})=\pi_2^*(c_{r+1})-2
 \pi_2^*(c_{r-1})\eta \theta + \bigl((rs+r)\eta+ \gamma \bigr)\pi_2^*(c_r)
\end{equation}
 and
 \begin{equation}\label{c4d2}
c_{r+2}(\pi_2^*(\cM)^{\vee}-\cB^{\vee})=\pi_2^*(c_{r+1})\bigl((rs+r)
\eta+ \gamma\bigr)-2\pi_2^*(c_r) \eta \theta.
\end{equation}
\end{remark}

\begin{proposition}\label{R1}
Let $[C]\in \cM_{rs+s}$ be a Brill-Noether general curve and denote
by $\P$ the Poincar\'e bundle on $C\times \mathrm{Pic}^{d}(C)$. We
have the following identities in $H^*(\rm{Pic}$$^{d}(C), \mathbb
Z)$: $$ c_1\big(R^1\pi_{2
*}(\P_{| C\times W^r_{d}(C)})\bigr)=\theta-c_1 \ \mbox{ and }\  c_2\bigl(R^1\pi_{2
*}(\P_{| C \times W^r_{d}(C)})\bigr)=\frac{\theta^2}{2}-\theta c_1+ c_2.$$
\end{proposition}
\begin{proof}
We recall that in order to obtain a determinantal structure on
$W^r_{d}(C)$ one fixes a divisor $D\in C_e$ of sufficiently high degree $e>0$ and
considers the morphism $$(\pi_2)_*\bigl(\P\otimes
\OO(\pi_1^*D)\bigr)\rightarrow (\pi_2)_*\bigl(\P\otimes
\OO(\pi_1^*D_{ |\pi_1^*D})\bigr).$$ Then $W^r_{d}(C)$ is the
degeneration locus of rank $d-g-r+e$ of this map and there is an
exact sequence of vector bundles over $W^r_{d}(C)$:
$$0\rightarrow \cM\longrightarrow (\pi_2)_*\bigl(\P\otimes
\OO({\pi_1}^*D)\bigr)\longrightarrow (\pi_2)_*\bigl(\P\otimes
\OO(\pi_1^*D)_{| \pi_1^*D}\bigr)\rightarrow R^1\pi_{2
*}\bigl(\P_{| C\times W^r_{d}(C)}\bigr)\longrightarrow 0.$$
\noindent
From this sequence our claim follows if we take into account that
$(\pi_2)_*\bigl(\P\otimes \OO(\pi_1^*D)_{| \pi_1^*D}\bigr)$ is
numerically trivial and $c_t\bigl((\pi_2)_*(\P \otimes
\OO(\pi_1^*D))\bigr)=e^{-\theta}.$
\end{proof}
\begin{remark}
For future reference we note that Proposition \ref{R1} provides a
quick way to compute the canonical class $K_{W^r_d(C)}$. Indeed,
since $T_{\mathrm{Pic}^d(C)}$ is trivial, we have that
$K_{W^r_d(C)}=c_1(N_{W^r_d(C)/\mathrm{Pic}^d(C)})$. From the
realization of $W^r_d(C)$ as a determinantal variety, we obtain that
$N_{W^r_d(C)}=Hom\bigl(\cM, R^1\pi_{2 *}\bigl(\P_{| C\times
W^r_d(C)}\bigr)\bigr)$, which leads to the expression:
\begin{equation}\label{canonical}
K_{W^r_d(C)}\equiv (r+1)\theta+(s-r-2)c_1.
\end{equation}

We shall also need in Section 3 the expressions for $K_X$ and $K_Y$.
To start with the surface $X$, we have that $K_X\equiv
(2rs+2s-2)\eta+K_{W^r_d(C)}+c_1(N_{X/C\times W^r_d(C)})$. Next we
use Proposition \ref{xy}, to express the normal bundle of the
determinantal subvariety $X\subset C\times W^r_d(C)$ as
$N_{X/C\times W^r_d(C)}=Hom(\mbox{Ker}(\zeta), \mbox{Coker}(\zeta))$
which leads to the formula:
\begin{equation}\label{canonicalx}
K_X\equiv (r+1)\theta+
(r-1)c_1(\mbox{Ker}(\zeta)^{\vee})+(s-r-1)\pi_2*(c_1) +2\gamma+
(6rs+2r+4s-2)\eta.
\end{equation}
In a similar manner, using the vector bundle map $\chi$, we find the
canonical class of $Y$:
\begin{equation}\label{canonicaly}
K_Y\equiv
(r+1)\theta+(r-1)c_1(\mbox{Ker}(\chi)^{\vee})+(s-r-1)\pi_2^*(c_1)
+\gamma+ (3rs+r+2s-2)\eta.
\end{equation}

\end{remark}

As a first step towards computing $[\gp_{g, d}^r]$
we determine the $\delta_1$ coefficient in its expression. For simplicity we set $$\tilde{\mathbb E}:=\mathbb E\otimes \OO_{\mm_g}\bigr(\sum_{j=1}^{[g/2]} (2j-1)\cdot \delta_j\bigr)$$ for the twist of the Hodge bundle.
\begin{theorem}\label{d1}
Let $[C]\in \cM_{rs+s}$ be a Brill-Noether general curve and denote
by $C^1\subset \Delta_1$ the associated test curve. Then  the
coefficient of $\delta_1$ in the expansion of $\gp_{g, d}^r$ in
terms of the generators of $\mathrm{Pic}(\mm_g)$ is  equal to
$$b_1=
\frac{rs(r+1)(s-1)\
C_{r+1}}{2(s+r+1)(s+r)(s+r+2)(rs+s-1)}\Bigl((2s^2+2s^3)r^4+(2s^4+12s^3+23s^2+9s)r^3+$$
+$$(8s^4+39s^3+75s^2+46s+10)r^2+
(10s^4+59s^3+108s^2+89s+26)r+4s^4+30s^3+64s^2+58s+12)\Bigr).$$
\end{theorem}
\begin{proof} We intersect the degeneracy locus of the map
$\F\otimes \cN\rightarrow \sigma^*(\tilde{\mathbb E})$ with the surface
$\sigma^*(C^1)$ and use that the vector bundles $\F$ and $\cN$ were
defined by retaining the sections of the genus $g-1$ aspect of each
limit linear series and dropping the information coming from the
elliptic curve. It follows that $Z_i\cdot c_2(\sigma^*(
\tilde{\mathbb{E}})-\F\otimes \cN)=0$ for $1\leq i\leq 3$ because both
$\sigma^*{\tilde{\mathbb E}}$ and $\F\otimes \cN$ are trivial along the
surfaces $Z_i$.  Furthermore, we also have that $[X_1\times
X_2]\cdot c_2(\sigma^*(\tilde{\mathbb E})-\F\otimes \cN)=0$, because
$c_2(\sigma^*(\tilde{\mathbb E})-\F\otimes \cN)_{| X_1\times X_2}$ is in
fact the pull-back of a codimension $2$ class from the
$1$-dimensional cycle $X_1$, therefore the intersection number is
$0$ for dimensional reasons. We are left with estimating the
contribution coming from $X$ and we write
$$\sigma^*(C^1)\cdot
c_2(\sigma^* \tilde{\mathbb E}-\F\otimes \cN)=c_2(\sigma^*\tilde{\mathbb E}_{|X})-c_1(\sigma^*\tilde{\mathbb E}_{|X})\cdot c_1(\F\otimes \cN_{|
X})+c_1^2(\F\otimes \cN_{| X})-c_2(\F\otimes \cN_{|X})$$ and we are
going to compute each term in the right-hand-side of this
expression.

Since we have a canonical identification $\tilde{\mathbb E}_{|
C^1}[C_y^1]=H^0(C, K_C\otimes \OO_C(2y))$ for each $y\in C$, we
obtain that $c_2(\sigma^*\tilde{\mathbb E}_{|X})=0$ and $c_1(\sigma^*
\tilde{\mathbb E}_{|X})=-(2g-4)\eta$. Recall also that we have set
$c_i(\F_{|X}^{\vee})=\pi_2^*(c_i)$ for $0\leq i\leq r+1$, where $c_i\in H^{2i}(W^r_d(C), \mathbb Z)$.

In Proposition \ref{xy} we introduced a vector bundle morphism
$\zeta: J_1(\mathcal{P})^{\vee}\rightarrow \pi_2^*(\cM)$ over
$C\times W^r_d(C)$. We denote by $U:=\mbox{Ker}(\zeta)$ and we view
$U$ as a line bundle over $X$ with fibre over a point $(y, L)\in X$
being the space
$$U(y, L)=\frac{H^1(C, L\otimes \OO_C(-2y))^{\vee}}{H^1(C,
L)^{\vee}}.$$ The Chern numbers of $U^{\vee}$ can be computed from
the Harris-Tu formula and we find that for any class $\xi \in
H^2(C\times W^r_d(C))$ we have the following (cf. (\ref{c3d1})):
$$c_1(U^{\vee})\cdot \xi_{|
X}=c_{r+1}\bigl(\pi_2^*(\cM)^{\vee}-J_1(\P)^{\vee}\bigr)\cdot
\xi_{|X}=$$ $$= \bigl(\pi_2^*(c_{r+1})-6\eta \theta \pi_2^*(c_{r-1})
+((4rs+2r+2s)\eta+2\gamma)\pi_2^*(c_r)\bigr)\cdot \xi_{|X},$$ and
$$c_1^2(U^{\vee})=c_{r+2}\bigl(\pi_2^*(\cM)^{\vee}-J_1(\P)^{\vee}\bigr)=\pi_2^*(c_{r+1})((4rs+2r+2s)\eta+2\gamma)-
6\pi_2^*(c_r)\eta \theta.$$ The line bundle $U$ is used to evaluate
the Chern numbers of $\cN_{|X}$ via the exact sequence:
\begin{equation}\label{ne}
0\longrightarrow \pi_2^*R^1\pi_{2 *}\bigl(\mathcal{P}_{|C\times
W^r_d(C)}\bigr)^{\vee} \longrightarrow \cN_{| X}\longrightarrow
U\longrightarrow 0,
\end{equation}
from which we obtain (by also using Proposition \ref{R1}), that
$c_1(\cN_{| X})=-\theta+c_1-c_1(U^{\vee})$ and $$c_2(\cN_{|
X})=c_2-\theta\cdot c_1+\frac{\theta^2}{2}+(\theta-c_1)\cdot
c_1(U^{\vee}).$$ Therefore we can write that
$$\sigma^*(C^1)\cdot c_2(\sigma^*\tilde{\mathbb E}-\F\otimes \cN)=(2g-4)\eta
\cdot c_1(\F\otimes \cN_{|X})+c_1^2(\F\otimes \cN_{|X})-c_2(\F
\otimes \cN_{|X})=$$
$$={r+2\choose 2}c_1^2(U^{\vee})+\Bigl((r+1)^2\cdot
\theta+2(r+1)(1-s(r+1))\cdot \eta +((r+1)(s-r-1)+1)\cdot
c_1\Bigr)c_1(U^{\vee})+\cdots,$$ where the term we omitted is a
quadratic polynomial in $\theta, \eta$ and $\gamma$ which will be
multiplied by the class $[X]$. Since we have already computed
$c_1(U^{\vee})$ and $c_1^2(U^{\vee})$, we can write
$\sigma^*(C^1)\cdot c_2(\sigma^* (\tilde{\mathbb E})-\F\otimes \cN)$ as a
polynomial in the classes $\pi^*_2(c_i)$, $\eta$, $\theta$ and
$\gamma$ and the only non-zero terms will be those which contain
$\eta$. Then we apply Propositions \ref{intersection} and
\ref{intersection2} and finally compute the coefficient
$$b_1:=\sigma^*(C^1)\cdot c_2(\sigma^*(\tilde{\mathbb E})-\F\otimes
\cN)/(2g-4),$$ which finishes the proof.

\end{proof}

\begin{theorem}\label{d0}
Let $[C, q]\in \cM_{rs+s, 1}$ be a Brill-Noether general pointed
curve and denote by $C^0\subset \Delta_0$ the associated test curve.
Then the $\delta_0$-coefficient of $[\gp_{g, d}^r]$ is given by the
formula:
$$b_0=\frac{r(r+2)(s-1)(s+1)(2rs+2s+1)(rs+s+2)(rs+s+6)}{12(rs+s-1)(s+r+2)(s+r)}C_{r+1}.$$
\end{theorem}

\begin{proof} We look at the virtual degeneracy locus
of the morphism $\F\otimes \cN \rightarrow \sigma^*(\tilde{\mathbb E})$
along the surface $\sigma^*(C^0)$. The first thing to note is that
the vector bundles $\F_{| \sigma^*(C^0)}$ and $\cN_{|
\sigma^*(C^0)}$ are both pull-backs of vector bundles on $Y$. For
convenience we denote this vector bundles also by $\F$ and $\cN$,
hence to use the notation of Proposition \ref{limitlin0}, $\F_{|
\sigma^*(C^0)}=\epsilon^*(\F_{| Y})$ and $\cN_{|
\sigma^*(C^0)}=\epsilon^*(\cN_{| Y})$. We find that
$$\sigma^*(C^0)\cdot c_2(\sigma^*\tilde{\mathbb E}-\F\otimes \cN)=c_2(\sigma^*\tilde{\mathbb E}_{|Y})-c_1(\sigma^*\tilde{\mathbb E}_{| Y})\cdot
c_1(\F\otimes \cN_{| Y})+c_1^2(\F\otimes \cN_{| Y})-c_2(\F\otimes
\cN_{|Y}),$$ and as in the proof of Theorem \ref{d1}, we are going
to compute each term in this expression. We denote by
$V:=\mbox{Ker}(\chi)$, where $\chi: \cB^{\vee}\rightarrow
\pi_2^*(\cM)^{\vee}$ is the bundle
 morphism coming from Proposition \ref{xy}. Thus $V$ is a line bundle on $Y$ with  fibre
$$V(y, L)=\frac{H^1(C, L\otimes \OO_C(-y-q))^{\vee}}{H^1(C, L)^{\vee}}$$
over each point $(y, L)\in Y$. By using again the Harris-Tu Theorem,
we find the following formulas for the Chern numbers of $V^{\vee}$
(cf. (\ref{c3d2}) and (\ref{c4d2})):
$$c_1(V^{\vee})\cdot \xi_{|Y}=c_{r+1}\bigl(\pi_2^*(\cM)^{\vee}-\cB^{\vee}\bigr)\cdot \xi_{| Y}=$$
$$=\bigl(\pi_2^*(c_{r+1})+\pi_2^*(c_r)((d-1)\eta+\gamma)-2\pi_2^*(c_{r-1})\eta
\theta\bigr)\cdot \xi_{| Y},$$ for any class $\xi\in H^2(C\times
W^r_{d}(C))$, and
$$c_1^2(V^{\vee})=c_{r+2}\bigl(\pi_2^*(\cM)^{\vee}-\cB^{\vee}\bigr)=\pi_2^*(c_{r+1})((d-1)\eta+\gamma)-
2\pi_2^*(c_r)\eta \theta.$$ To evaluate the Chern numbers of
$\cN_{|Y}$ we fit the line bundle $V$ in the following exact
sequence:
\begin{equation}\label{ne1}
0\longrightarrow \pi_2^*R^1\pi_{2 *}\bigl(\mathcal{P}_{|C\times
W^r_d(C)}\bigr)^{\vee} \longrightarrow \cN_{| Y}\longrightarrow
V\longrightarrow 0.
\end{equation}
This allows us to compute $c_1(V^{\vee})$ and $c_1^2(V^{\vee})$ and
then we can write that
$$\sigma^*(C^0)\cdot c_2(\sigma^*\tilde{\mathbb E}-\F\otimes \cN)=\eta
\cdot c_1(\F\otimes \cN_{|Y})+c_1^2(\F\otimes \cN_{|Y})-c_2(\F
\otimes \cN_{|Y})=$$
$$={r+2\choose 2}c_1^2(V^{\vee})+\Bigl((r+1)^2\cdot
\theta+2(r+1)(r+1)\cdot \eta +((r+1)(s-r-1)+1)\cdot
c_1\Bigr)c_1(V^{\vee})+\cdots,$$ where the term we omitted is a
quadratic polynomial in $\theta, \eta$ and $\gamma$ which will be
multiplied by the class $[Y]$. Using repeatedly Propositions
\ref{intersection} and \ref{intersection2}, we finally evaluate all
the terms and obtain the stated expression for $b_0$ using the
relation $(2g-2)b_0-b_1= \sigma^*(C^0)\cdot c_2(\sigma^*\tilde{\mathbb
E}-\F\otimes \cN).$
\end{proof}

We finish the calculation of $s(\gp_{g, d}^r)$ by proving the
following result:

\begin{theorem}\label{elltail}
Let $[C, q]\in \cM_{rs+s, 1}$ be a suitably general pointed curve
and $L\in W^r_{d}(C)$ any linear series with a cusp at $q$. Then the
multiplication map $$\mu_0(L, q): H^0(C, L)\otimes H^0(K_C\otimes
L^{\vee}\otimes \OO_C(2q))\rightarrow H^0(C, K_C\otimes \OO_C(2q))$$
is injective. If $\gp_{g, d}^r\equiv a\lambda -\sum_{j=0}^{[g/2]}
b_j\delta_j \in \mathrm{Pic}(\mm_g)$, we have the relation
$a-12b_0+b_1=0$.
\end{theorem}
\begin{proof} We consider again the pencil $R\subset \mm_g$ obtained
by attaching to $C$ at the point $q$ a pencil of plane cubics. It is
well-known that $R\cdot \lambda=1, R\cdot \delta_0=12$ and $R\cdot
\delta_1=-1$, thus the relation $a-12b_0+b_1=0$ would be immediate
once we show that $R\cap \gp_{g, d}^r=\emptyset$. Assume by
contradiction that $R\cap \gp_{g, d}^r\neq \emptyset$ and then
according to Proposition \ref{limitlin1} there exists $L\in
W^r_{d}(C)$ with
 $h^0(L\otimes \OO_C(-2q))=r$ such that the multiplication map
$\mu_0(L, q)$ is not injective.

We degenerate $[C, q]$ to the stable curve $$[C_0:=E_0\cup_{p_1} E_1\cup _{p_2}
E_2\cup \ldots \cup_{p_{g-2}} E_{g-2}, p_0]\in \mm_{g-1, 1}$$
consisting of a string of elliptic curves such that $p_0\in E_0$ and
the differences $p_{i+1}-p_i\in \mbox{Pic}^0(E_i)$ for $0\leq i\leq
g-3=rs+s-2$, are not torsion classes.
 For each $0\leq i \leq g-2$ we
denote by $L^i\in \mbox{Pic}^d(C_0)$ the unique limit of the line
bundles  $L_t\in \mbox{Pic}^d(C_t)$ having the property that
$\mbox{deg}(L^i_{| E_j})=0$ for $i\neq j$. Here $[C_t, p_t]\in
\cM_{g-1, 1}$ is a $1$-dimensional family of smooth pointed curves
with the property $\mbox{lim}_{t\rightarrow 0} [C_t, p_t]=[C_0,
p_0]\in \mm_{g-1, 1}$ and where we also assume that $\mbox{Ker}\
\mu_0(L_t, p_t)\neq 0$ for all $t\neq 0$.

 Similarly, we define
 $M^i \in \mbox{Pic}^{2g-2-d}(C_0)$ to be the unique limit of the line bundles
 $K_{C_t}\otimes L_t^{\vee}\otimes
\OO_{C_t}(2p_t)$  characterized by the property $\mbox{deg}(M^i_{|
E_j})=0$ for $i\neq j$. We denote by $\{(L^i_{| E_i},
V_i)\}_{i=0}^{g-2}$ and by $\{(M^i_{|E_i}, W_i)\}_{i=0}^{g-2}$ the
limit linear series on $C_0$ corresponding to $L_t$ and
$K_{C_t}\otimes L_t^{\vee}\otimes \OO_{C_t}(2p_t)$ respectively as
$t\rightarrow 0$. Reasoning along the lines of \cite{EH3} or
\cite{F1} Proposition 3.2, for each $0\leq i\leq g-2$ we find
elements
$$0\neq \rho_i\in \mbox{Ker}\{V_i\otimes W_i\rightarrow H^0(E_i,
L^i\otimes M^i_{| E_i})\}$$ satisfying
$\mbox{ord}_{p_{i+1}}(\rho_{i+1})\geq
\mbox{ord}_{p_{i}}(\rho_{i})+2$ for $0\leq i\leq g-3$. Since
$\mbox{ord}_{p_0}(\rho)\geq 2$, we find that
$\mbox{ord}_{p_{g-2}}(\rho_{g-2})\geq 2g-2=\mbox{deg}(L^{g-2}\otimes
M^{g-2}_{| E_{g-2}})$, which is impossible.
\end{proof}

\section{Maps between moduli spaces of curves}

We begin the study of the map $\phi: \mm_g-->\mm_{g'}$ given by
$\phi([C]):=[W^r_d(C)]$ in the case $\rho(g, r, d)=1$, so that
$g=rs+s+1$ and $d=rs+r+1$. The genus of $W^r_d(C)$ for a general
$[C]\in \cM_g$ has been computed in \cite{EH2} Theorem 4, and we
have the formula:
\begin{equation}\label{genusbig}
g'=g(W^r_d(C))=1+g!\ \frac{s(r+1)}{s+r+1}\prod_{i=0}^r
\frac{i!}{(s+i)!}.
\end{equation}
We shall describe the pull-back map $\phi^*:
\mbox{Pic}(\mm_{g'})\rightarrow \mbox{Pic}(\mm_g)$ and to avoid
confusion we denote, as usual, by $\lambda, \delta_0, \ldots,
\delta_{[g/2]}$ the generators of $\mbox{Pic}(\mm_g)$, and by
$\lambda', \delta_0', \ldots, \delta'_{[g'/2]}$ the generators of
$\mbox{Pic}(\mm_{g'})$. We start by describing the map $\phi$ over a
generic point of each boundary divisors $\Delta_j$ for $0\leq j\leq
[g/2]$. If $[C^j_y:=C\cup_y D]\in \Delta_j$ is a general point, then
$\phi([C^j_y])$ is the stable reduction of the variety
$\overline{G}^r_d(C^j_y)$ of limit linear series $\mathfrak g^r_d$.
Our analysis shows that $\overline{G}^r_d(C^j_y)$ is always a
semi-stable curve and this observation completely determines $\phi$
in codimension $1$.

Suppose that $[C]\in \cM_{rs+s}$ is a Brill-Noether-Petri  general
curve and that $[E, y]\in \cM_{1, 1}$ is a pointed elliptic curve.
We recall that we have introduced the smooth surface $X=\{(y, L)\in
C\times W^r_d(C): h^0(C, L\otimes \OO_C(-2y))=r\}$. For $y\in C$ we
denote by $X_y:=\pi^{-1}_1(y)$ the fibre of the first projection
$\pi_1:X\rightarrow C$. For each of the $n_0$ linear series $A\in
W^r_{d-1}(C)$ there exists a section $\sigma_A:C\rightarrow X$ given
by $\sigma_A(y)=(y, A\otimes \OO_C(y))$ and we set
$\Sigma_A:=\mbox{Im}(\sigma_A)$. From the description given in
Proposition \ref{limitlin1}, it follows that $\phi([C\cup_y E])$ is
the stable curve of genus $g'$ obtained by attaching to the spine
$X_y$ copies of $E\cong \mbox{Pic}^{r+1}(E)$ at the points
$\sigma_A(y)$ for each $A\in W^r_{d-1}(C)$.

Similarly, having fixed a general pointed curve $[C, q]\in
\cM_{rs+s, 1}$, we recall that we have introduced the surface
$Y=\{(y, L)\in C\times W^r_d(C): h^0(L\otimes \OO_C(-y-q))=r\}$.
(cf. Proposition \ref{limitlin0}). For each $y\in Y$ we set
$Y_y:=\pi_1^{-1}(y)$. To every $A\in W^r_{d-1}(C)$ correspond two
sections $u_A: C\rightarrow Y$, $v_A:C \rightarrow Y$ given by
$u_A(y)=(y, A\otimes \OO_C(y))$ and $v_A(y)=(y, A\otimes \OO_C(q))$
respectively. If as before, we denote by $[C^0_y]:=[C/y\sim q]\in
\Delta_0$, then $\phi([C_y^0])$ is the stable curve obtained from
$Y_y$ by identifying the points $v_A(y)$ and $u_A(y)$ for all linear
series $A\in W^r_{d-1}(C)$.

For $2\leq j\leq [g/2]$, the irreducible components of
$\phi([C^j_y])$ are indexed by Schubert indices
$\overline{\alpha}:=(\alpha_0\leq \ldots \leq \alpha_r)$ such that
there exists limit linear series $l=\{l_C, l_D\}\in
\overline{G}^r_d(C^j_y)$ with $\alpha^{l_C}(y)=\overline{\alpha}$,
$\rho(l_C, y)\in \{0, 1\}$ and $\rho(l_C, y)+\rho(l_D, y)=1$ (a
precise list of such $\overline{\alpha}$'s is given in the proof of
Theorem \ref{divisorspullback}). To describe the pull-backs of the
tautological classes under $\phi$ we need a description of the
numerical properties of the push-forwards under $\phi$ of the
standard test curves $R$ and $C^j$ where $0\leq j\leq [g/2]$. We
carry this out in detail only for $j=0, 1$ which is sufficient to
compute the slopes of pull-backs $\phi^*(D')$ where $[D']\in
\mbox{Pic}(\mm_{g'})$. The case $j\geq 2$ is quite similar and again
we skip these details. To keep our formulas relatively simple we
only deal with the case $r=1$, when $g=2s+1$ and
$$\phi:\mm_{2s+1}-->\mm_{1+\frac{s}{s+1}{2s+2\choose s}}.$$

\begin{proposition}\label{pushforwardr}
We fix a general pointed curve $[C, q]\in \cM_{2s, 1}$ and we
consider the test curve $R\subset \mm_{2s+1}$ obtained by attaching
a pencil of plane cubics to $C$ at the fixed point $q$. If $n_0:=\#(W^1_{s+2}(C))$,
then we have the following relations:
$$\phi_*(R)\cdot \lambda'=n_0,$$
$$ \phi_*(R)\cdot \delta'_0=12n_0, $$
$$ \phi_*(R)\cdot \delta'_1=-n_0,$$
and
$$\phi_*(R)\cdot \delta'_j=0 \mbox{ for } j\geq 2.$$
\end{proposition}
\begin{proof} We denote by $f:\widetilde{\PP}^2:=\mbox{Bl}_9(\PP^2)\rightarrow \PP^1$ the
fibration induced by a pencil of plane cubics after blowing-up the
$9$ base points of the pencil. Since $f$ has $9$ sections, there is
an isomorphism between $f$ and its Picard fibration
$\mbox{Pic}^2(f)\rightarrow \PP^1$. The curve $\phi_*(R)\subset
\mm_{g'}$ is induced by a fibration of stable curves  $\pi:
T\rightarrow \PP^1$, where $$\pi^{-1}(t)=X_q\bigcup \bigl\{\bigcup
_{\sigma_A(q)} \mbox{Pic}^2(f^{-1}(t)):A\in W^r_{d-1}(C)\bigr\}, \mbox{
 } \mbox{ for each } t\in \PP^1.$$ In other words, $\pi$ is obtained by
attaching to the fixed curve $X_q$, $n_0$ copies of the elliptic
curve $f^{-1}(t)$ at each of the points $\sigma_A(q)$. The claimed
formulas are now immediate.
\end{proof}

\begin{proposition}\label{pushforward1} We fix general curves
$[C]\in \cM_{2s}$ and $[E, y]\in \cM_{1, 1}$ and consider the
associated test curve $C^1\subset \Delta_1\subset \mm_{2s+1}$. Then
we have the formulas
$$\phi_*(C^1)\cdot \lambda'=n_0 \frac{2s(s-1)(6s^2+10s+1)}{s+2},$$
$$ \phi_*(C^1)\cdot \delta'_0=C^1\cdot
\gp_{2s+1, s+2}^1, \mbox{  } $$ $$ \phi_*(C^1)\cdot
\delta'_1=-n_0(4s-2),$$
$$ \phi_*(C^1)\cdot \delta'_j=0
\mbox{ for } j\geq 2.$$
\end{proposition}
\begin{proof}
The $1$-cycle $\phi_*(C^1)$ corresponds to a family of curves
constructed as follows: We start with $\pi_1:X\rightarrow C$ and
consider the sections $\{\sigma_A:C\rightarrow X\}_{A\in
W^1_{s+1}(C)}$. We also consider $n_0$ disjoint copies of the
trivial family $C\times E\rightarrow C$ which we glue to $\pi_1$
along each of the sections $\sigma_A$. From this description it
follows that $\phi_*(C^1)\cdot \delta_0' =C^1\cdot \gp_{2s+1,
s+2}^1$ and this equals the number of points $y\in C$ (counted with
the appropriate multiplicities) such that $X_y$ is singular at some
point $(y, L)$. This  translates into saying that the Petri map
$\mu_0(L, y): H^0(C, L)\otimes H^0(K_C\otimes L^{\vee}\otimes
\OO_C(2y))\rightarrow H^0(C, K_C\otimes \OO_C(2y))$ is not
injective. Also, $\phi_*(C^1)\cdot \delta'_j=0$ for $j\geq 2$. Using
the description of $N_{\Delta'_1/\mm_{g'}}$ we have that
$$\phi_*(C^1)\cdot \delta_1'=\sum_{A\in W^1_{s+1}(C)}
(\Sigma_A)^2=\sum_{A\in W^1_{s+1}(C)} \Bigl(2g(C)-2-\Sigma_A\cdot
K_X\Bigr).$$ To estimate this sum, we recall that we have computed
the canonical class of $X$ (cf. (\ref{canonicalx})):
$$K_X\equiv 2\theta+(s-2)\cdot \pi_2^*(c_1)+2\gamma+10s\cdot \eta. $$ By direct
computation we obtain that $\Sigma_A\cdot
\theta=\sigma_A^*(\theta)=2s, \Sigma_A\cdot \eta=1$ and
$\Sigma_A\cdot \gamma=-4s$ (all these intersection numbers are being
computed on the smooth surface $X$).

We now compute $\Sigma_A\cdot \pi_2^*(c_1)$ and note that
$\sigma_A^*\pi_2^*(\cM)=(p_2)_*\bigl(\mu^*(\P_{| C\times
W^r_d(C)})\bigr)$, where $\mu:C\times C\rightarrow C\times
\mbox{Pic}^d(C)$ is defined as $\mu(x, y)=(x, A\otimes \OO_C(y))$
and $p_1, p_2:C\times C\rightarrow C$ are the two projections. The
key observation here is that if the Poincar\'e bundle $\P$ is chosen
in such a way that $\P_{| \{q\}\times \mathrm{Pic}^d(C)}$ is trivial
for a point $q\in C$, then $$\mu^*(\P)=p_1^*(A)\otimes \OO_{C\times
C}(\Delta)\otimes p_2^*(\OO_C(-q)),$$ hence
$(p_2)_*\mu^*(\P)=(p_2)_*\bigl(p_1^*A\otimes \OO_{C\times
C}(\Delta)\bigr)\otimes \OO_C(-q)$. Then we note that  the vector
bundle $(p_2)_*\bigl(p_1^*(A)\otimes \OO_{C\times C}(\Delta)\bigr)$
is trivial, thus it  follows that $\mbox{deg}
\bigl(\sigma_A^*\pi_2^*(\cM^{\vee})\bigr)=\mbox{rank}(\cM)=h^0(A)=2$.
Putting these calculations together, we obtain that $\Sigma_A \cdot
K_X=8s-4$ and then $(\Sigma_A^2)=-4s+2$, that is, $\phi_*(C^1)\cdot
\delta_1'=-n_0(4s-2)$.

We are left with the computation of $\phi_*(C^1)\cdot \lambda'$,
which equals the degree of the Hodge bundle over the family
$\pi_1:X\rightarrow C$. From the Mumford relation
$\kappa_1=12\lambda-\delta$  we find that
$$\phi_*(C^1)\cdot
\lambda'=\frac{K_{X/C}^2+\delta(\pi_1)}{12}=\frac{K_{X/C}^2+C^1\cdot
\gp_{2s+1, s+2}^1}{12},$$ where $K_{X/B}=K_X-\pi_1^*(K_C)$ is the
relative canonical class. A direct calculation involving
Propositions \ref{intersection} and \ref{intersection2} shows that
$$K_{X/C}^2=(6s^3-10s^2-4s)\pi_2^*(c_1^2)+(24s^2-32s+16)
\pi_2^*(c_1)\cdot \theta+24s^2\theta^2.$$ The calculation of
$[\gp_{2s+1, s+2}^1]$ (precisely Theorem \ref{d1}), yields that
$C^1\cdot \gp_{2s+1, s+2}^1=4n_0s(s-1)(12s^2+23s+8)/(s+2)$, which
leads to the stated formula for $\phi_*(C^1)\cdot \lambda'$.
\end{proof}

\begin{proposition}\label{pushforward0}
We fix a general pointed curve $[C, q]\in \cM_{2s, 1}$ and consider
the test curve $C^0\subset \Delta_0\subset \mm_{2s+1}$. Then we have
the following formulas:
$$\phi_*(C^0)\cdot \delta_1'=n_0,$$ 
$$\phi_*(C^0)\cdot
\delta_0'=C^0\cdot \gp_{2s+1, s+2}^1- 4n_0s,$$
$$\phi_*(C^0)\cdot \lambda'=n_0\frac{s(s-1)(2s^2+4s+1)}{s+2}$$
$$\phi_*(C^0) \cdot
\delta'_j=0 \ \mbox{ for } j\geq 2.$$
\end{proposition}
\begin{proof} We describe the family of stable curves
inducing $\phi_*(C^0)$. We start with the family $\pi_1:Y\rightarrow
C$ and consider the sections $\{u_A, v_A:C\rightarrow Y\}_{A\in
W^1_{s+1}(C)}$ with images $U_A:=u_A(C)$ and $V_A:=v_A(C)$
respectively. We denote by $Y'$ the blow-up of $Y$ at the $n_0$
points of intersections $\{U_A\cap V_A=(q, A\otimes
\OO_C(q))\}_{A\in W^1_{s+1}(C)}$ (see also Proposition
\ref{limitlin0}), and we denote by $\widetilde{U}_A$ and
$\widetilde{V}_A$ the strict transforms of $U_A$ and $V_A$
respectively. Then $\phi_*(C^1)\subset \mm_{g'}$ is induced by the
fibration $\pi:\widetilde{Y}\rightarrow C$, where
$\widetilde{Y}=\sigma^*(C^0)$ is  the surface obtained from $Y'$ by
identifying the sections $\widetilde{U}_A$ and $\widetilde{V}_A$ for
each $A\in W^1_{s+1}(C)$. The numerical characters of $\phi_*(C^0)$
are now easily describable. We have that $\phi_*(C^0)\cdot
\delta'_j=0$ for $j\geq 2$,\  \ \ $\phi_*(C^0)\cdot \delta'_1=n_0$,
and
$$\phi_*(C^0)\cdot \delta'_0=C^0\cdot \gp_{2s+1, s+2}^1+\sum_{A\in
W^1_{s+1}(C)} \Bigl((U_A)^2+(V_A)^2-2\Bigr).$$ We recall that we
have computed the canonical class of $Y$ (cf. (\ref{canonicaly})):
$$K_Y\equiv 2\theta+(s-2)\pi_2^*(c_1)+\gamma+(5s-1)\cdot \eta.$$
Since $u_A^*(\theta)=g(C)=2s, \ u_A^*(\gamma)=-4s, \ u_A^*(\eta)=1$
and $u_A^*\pi_2^*(c_1)=2$ (the proof of this last equality follows
from the calculation in Proposition \ref{pushforward1}), we find
that $K_Y\cdot U_A=7s-5$, hence by the adjunction formula
$(U_A)^2=-3(s-1)$. Similarly, $(V_A)^2=-s-1$, therefore
$(U_A)^2+(V_A)^2-2=-4s$, for every $A\in W^1_{s+1}(C)$, which
determines $\phi_*(C^0)\cdot \delta_0'$.

We still have to estimate $\phi_*(C^0)\cdot \lambda'$. Like in
Proposition \ref{pushforward1}, using Mumford's relation, this
number equals the degree of the Hodge bundle on the family
$\pi_1:Y\rightarrow C$:
$$\phi_*(C^0)\cdot
\lambda'=\frac{K_{Y/C}^2+\delta(\pi_1)}{12}=\frac{K_{Y/C}^2+C^0\cdot
\gp_{2s+1, s+2}^1}{12}.$$ From Theorem \ref{d1} we know that
$C^0\cdot \gp_{2s+1, s+2}^1=2n_0s(s-1)(4s^2+9s+4)/(s+2)$. By direct
computation we also obtain that
$$K_{Y/C}^2=(s^3-s^2-2s)\pi_2^*(c_1^2)+(4s^2-4s+2)\pi_2^*(c_1)\cdot
\theta+4(s-1)\theta^2.$$ Moreover, $\pi_2^*(c_1^2)=n_0(4s+2)/(s+2),\
\pi_2^*(c_1)\cdot \theta =2n_0s$ \  and $\theta^2=n_0s(s+1)$ (all these
intersection numbers are being computed on $Y$ using Proposition
\ref{intersection}). This completes the calculation of
$\phi_*(C^0)\cdot \lambda'$.
\end{proof}

We are in a position to describe pull-backs of divisors classes
under the map $\phi$:

\begin{theorem}\label{divisorspullback}
 We consider the rational map $\phi:\mm_g --> \mm_{g'}$,
\  $\phi[C]=[W^1_{s+2}(C)]$, \ where
$$g:=2s+1 \mbox{  } \mbox{  and } \mbox{ } g':=1+\frac{s}{s+1}{2s+2\choose
s}.$$ We then have the following description of the map $\phi^*:\mathrm{Pic}(\mm_{g'})\rightarrow
\mathrm{Pic}(\mm_g)$:
$$\phi^*(\lambda')=n_0\Bigl(\frac{6s^4+20s^3-s^2-20s-2}{(s+2)(2s-1)}\lambda-\frac{s(s^2-1)}{2s-1}\delta_0-$$
$$-
\frac{2s(s-1)
(6s^2+10s+1)}{(s+2)(4s-2)}\delta_1-\sum_{i=2}^{[g'/2]} b_i \delta_i \Bigr), $$
where $$b_i\geq \frac{s(s^2-1)}{2s-1} \ \mbox{ for } \ 2\leq i\leq [g/2],$$
$$\phi^*(\delta_0')=n_0\cdot \delta_0+[\gp_{2s+1, s+2}^1], \
\phi^*(\delta_1')=n_0\cdot \delta_1 \ \mbox{ and } \phi^*(\delta_j')=0 \ \mbox{ for } 2\leq j\leq [g'/2].$$
\end{theorem}
\begin{proof}
The formulas involving $\phi^*(\lambda'), \phi^*(\delta_0')$ and
$\phi^*(\delta_1')$ are consequences of Propositions
\ref{pushforward1} and \ref{pushforward0} via the push-pull formula.
To prove that $\phi^*(\delta_j')=0$ for $2\leq j\leq [g'/2]$, we
show that $\mathrm{\phi}_*(\mm_g) \cap \Delta_j'=\emptyset$. This
follows once we note that (i) the generic point of every component
of $\gp_{2s+1, s+2}^1$ corresponds to a curve $[C]\in \cM_{2s+1}$
for which $W^1_{s+2}(C)$ is irreducible with precisely one node, and
that (ii) $\phi_*(\Delta_j)\subset \mm_{g'}- \bigcup_{i\geq
2}^{[g'/2]} \Delta_i'$ for every $2\leq j\leq [g/2]$. Indeed, let us
fix a general point $[C^j_y:=C\cup_y D]\in \Delta_j$ where $[C,
y]\in \cM_{j, 1}$ and $[D, y]\in \cM_{g-j, 1}$ are Brill-Noether
general pointed curves. For a real number $t$ we introduce the
notation $t_{+}:=\mathrm{max}\{t, 0\}$. The irreducible components
of the stable curve $\phi([C^j_y])$ are indexed by the set
$\mathcal{P}_j$ of Schubert indices
$$\overline{\alpha}:=\bigr(0\leq \alpha_0\leq \ldots  \leq \alpha_r\leq rs+1 \bigl) $$
satisfying the conditions (cf. \cite{EH2}, Proposition 1.2):
\begin{equation}\label{indexing}
\sum_{i=0}^r (\alpha_i+j-rs-1)\in \{j-1, j\}, \  \mbox{  }
\sum_{i=0}^r  (\alpha_i+j-rs-1)_{+}\leq j \mbox{ and } \sum_{i=0}^r
(g-j-\alpha_i)_{+}\leq g-j.
\end{equation}
For $\alpha\in \mathcal{P}_j$ we consider the (non-empty) variety
$\overline{G}^r_d(X)_{\overline{\alpha}}:=\{l\in
\overline{G}^r_d(X): \alpha^{l_C}(y)\geq \overline{\alpha} \}$ which
is a disjoint union of irreducible components of $\phi([X])$. When
$j\geq 2$, we claim that the stable curve $\phi([C^j_y])$ is not of
compact type. Using (\ref{indexing}) one checks that for every
$\alpha\in \mathcal{P}_j$ there are at least two partitions
$\overline{\beta}_1,\overline{\beta}_2\in \mathcal{P}_j$ such that
$\overline{G}^r_d(X)_{\overline{\alpha}} \cap
\overline{G}^r_d(X)_{\overline{\beta}_k}\neq \emptyset$ for $k=1,
2$. Thus for every component $Z$ of $\phi([C^j_y])$ we have that
$\#\bigl(Z\cap \overline{\phi([C^j_y])-Z} \bigr)\geq 2$, which
proves our claim.
\end{proof}

Theorem \ref{divisorspullback} contains enough information to encode
the slope of the pull-backs $\phi^*(D)$ for all classes $D\in
\mathrm{Pic}(\mm_{g'})$ and thus to prove Theorem \ref{slopepull}:
 If $s(D)=c$, then we have the following
formula for the slope of $\phi^*(D)\in \mathrm{Pic}(\mm_g)$:
$$s(\phi^*(D))=6+ \frac{8s^3(c-4)+5cs^2-30s^2+20s-8cs-2c+24}{s(s+2)(cs^2-4s^2-c-s+6)}\ .$$

\section{The map $\phi$ in small genus and applications to Prym varieties}

In this section we denote by $\cR_g$ the stack of \'etale double
covers of smooth curves of genus $g$ and by $\overline{\cR}_g$ its
compactification by means of Beauville admissible double covers, cf.
\cite{B}. It is proved in \cite{BCF} that $\overline{\cR}_g$ is
isomorphic to the stack parameterizing \emph{Prym curves} of genus
$g$, that is, data of the form $(X, L, \beta)$, where $X$ is a
quasi-stable curve with $p_a(X)=g$, $L \in \mathrm{Pic}^0(X)$ is a
line bundle such that $L_{| R}=\OO_R(1)$ for every destabilizing
rational component $R\subset X$ with $\#(R\cap (\overline{X-R}))=2$,
and $\beta: L^{\otimes 2}\rightarrow \OO_X$ is a sheaf homomorphism
whose restriction to the generic point of each component of $X$ is
non-zero. One has a finite  branched cover $\pi: \overline{\cR}_g
\rightarrow \mm_g$ and a regular morphism $\chi:
\overline{\cR}_g\rightarrow \mm_{2g-1}$ which assigns to an
admissible double cover the stable model of its source curve. We set
$\lambda:=\pi^*(\lambda)\in \mbox{Pic}(\overline{\cR}_g)$ and define
the following three irreducible boundary divisors in $
\overline{\cR}_g$:

\noindent $\bullet$ $\Delta_0^{'}$, with generic point being a Prym
curve $t:=[C_y^0:=C/y\sim q, L]$, where $[C, y, q]\in \cM_{g-1, 2}$
and $L \in \mathrm{Pic}^0(C_y^0)[2]$ is a line bundle such that if
$\nu:C_y^0\rightarrow C$ denotes the normalization map, then
$\nu^*(L)\neq \OO_C$. If $\tilde{C}\rightarrow C$ is the \'etale
$2:1$ cover induced by $\nu^*(L)$ and $y_i, q_i (i=1, 2)$ are the
inverse images of $y$ and $q$, then $\chi(t)=[\tilde{C}/y_1\sim q_1,
y_2\sim q_2]$.

\noindent $\bullet$ $\Delta_0^{''}$, with generic point
corresponding to $t:=[C_y^0, L]$ as above, but where
$\nu^*(L)=\OO_C$. In this case $\chi(t)$ consists of two copies
$[C_i, y_1, q_i]$ $(i=1, 2)$ of $[C, y, q]$, where we identify $y_1$
with $q_2$ and $y_2$ and $q_1$ respectively.

\noindent $\bullet$ $\Delta_0^r$, with generic point corresponding
to a Prym curve $t:=[X:=C\cup_{\{y, q\}} \PP^1, L]$, with $[C, y,
q]\in \cM_{g-1, 2}$ and  $L\in \mathrm{Pic}(C)$ is a line bundle
such that $L^{\otimes 2}=\OO_C(-y-q)$. In this case, if
$\tilde{C}\rightarrow C$ is the double cover induced by $L_{| C}$
and branched at $y$ and $q$ and if $\tilde{y}, \tilde{q}\in
\tilde{C}$ are the ramification points above $y$ and $q$
respectively, then $\chi(t)=[\tilde{C}/\tilde{y}\sim \tilde{q}]$.

For a straightforward dictionary between Beauville covers and Prym
curves we refer to \cite{D}. Note that
$\pi^*(\Delta_0)=\Delta_0'+\Delta_0^{''}+2\Delta_0^r$ and
$\Delta_0^r$ is the ramification locus of $\pi$. As usual, we set
$\delta_0^{'}:=[\Delta_0^{'}], \delta_0^{''}:=[\Delta_0^{''}]$ and
$\delta_0^r:=[\Delta_0^r]\in \mathrm{Pic}(\overline{\cR}_g)$. We
also denote by $p:\cC\rightarrow \overline{\cR}_g$ the universal
curve and by $\L$ the line bundle over $\cC$ whose restriction to
each fibre of $p$ is the underlying line bundle corresponding to a
Prym curve. In \cite{F3}, for each $i\geq 1$  we introduce the
tautological vector bundles $\mathbb E_i:=p_*(\omega_p\otimes
\L^{\otimes i})$ over $\overline{\cR}_g$ and we show that
\begin{equation}\label{prymtauto}
c_1(\mathbb E_i)={i\choose 2} \pi^*(\kappa_1)+\lambda-\frac{i^2}{4}
\delta_0^r.
\end{equation}

We discuss the geometry of the rational map $\phi:\mm_g-->\mm_{g'}$
for small values of $g=2s+1$. When $s=1$, then $g=g'=3$ and the map
$\phi:\mm_3\rightarrow \mm_3$ is simply the identity. Indeed, for a
smooth curve $[C]\in \cM_3$, we have a natural isomorphism $C\cong
W^1_3(C)$ given by $C\ni y\mapsto K_C\otimes \OO_C(-y)$ (Note that
this isomorphism extends over the hyperelliptic locus as well, when
$W^1_3(C)=C+W^1_2(C)$).

The first truly interesting case is $s=2$, when we have a map
$\phi:\mm_5-->\mm_{11}$
$$\phi([C]):=[W^1_4(C)]=[\mathrm{Sing}(\Theta_C)].$$ By
duality there is an involution $\tau: W^1_4(C)\rightarrow W^1_4(C)$
given by $\tau(L)=K_C\otimes L^{\vee}$. For $[C]\in
\cM_5-\mathcal{GP}_{5, 4}^{1}$ (that is, when $C$ is not trigonal
and possesses no vanishing theta-nulls), $\tau$ has no fixed points
and it induces an \'etale $2:1$ cover $f:W^1_4(C)\rightarrow
\Gamma$, where $[\Gamma]\in \cM_6$. Therefore $\phi$ factors to give
a map $\nu:\mm_5 --> \overline{\mathcal{R}}_6$. Moreover, there is
an isomorphism of principally polarized abelian varieties of
dimension $5$:
$$\bigl(\mathrm{Prym}(W^1_4(C)/\Gamma), \Xi\bigr)\cong \bigl(\mathrm{Jac}(C),
\Theta_C\bigr)$$ (see \cite{ACGH} pg. 296-301 or \cite{DS} for
details on this classically understood situation). The genus $6$
curve $\Gamma$ is identified with the locus of rank $4$ quadrics
containing the canonical curve $C\subset \PP^4$, and if $Q\in
\Gamma$ is such a quadric, then $f^{-1}(Q)$ consists of the
$\mathfrak g^1_4$'s determined by the two rulings on $Q$. If $[C]\in
\cM_5-\mathcal{GP}_{5, 4}^{1}$ then $\Gamma$ is a smooth plane
quintic. When $[C]\in \mathcal{GP}_{5, 4}^{1, 0}$, the curve
$\Gamma$ has nodes at the points corresponding to quadrics of rank
$3$. We have the following result which completely determines $\phi$
in codimension $1$:
\begin{proposition}\label{quintic}
The image of the rational map $\phi:\mm_5-->\mm_{11}$ given by
$\phi([C])=[W^1_4(C)]$ equals the closure
$\overline{\cM\mathcal{Q}}^{+}$ of the locus of genus $11$ curves
which are even \'etale double covers of smooth plane quintic curves.

\noindent $\bullet$ For a trigonal curve $[C]\in \cM_{5, 3}^1$, if
$A\in W^1_3(C)$ denotes the unique $\mathfrak g^1_3$, then
$\phi([C])$ consists of two copies of $C$ joined together at two
points $x, y\in C$ such that $x+y= |K_C\otimes A^{\otimes (-2)}|$.

\noindent $\bullet$ For a curve $[C\cup_y E]\in \Delta_1\subset
\mm_5$ where $g(C)=4$ and $g(E)=1$, $\phi([C\cup_y E])$ is a stable
curve of compact type consisting of a genus $9$ spine $\{L\in
W^1_4(C): h^0(L\otimes \OO_C(-2y))\geq 1\}$ and two elliptic tails
isomorphic to $E$ attached at the points $A\otimes \OO_C(y)$ where
$A\in W^1_3(C)$.

\noindent $\bullet$ For a curve $[C/y\sim q]\in \Delta_0\subset \mm_5$ where $[C, y, q]\in \cM_{4, 2}$,
$\phi([C/y\sim q])$ is the irreducible stable curve obtained from the smooth genus $9$ curve
$\{L\in W^1_4(C): h^0(L\otimes \OO_C(-y-q))\geq 1\}$ by identifying the two pairs of points
$A\otimes \OO_C(y)$ and $A\otimes \OO_C(q)$ for every $A\in W^1_3(C)$.

\noindent $\bullet$ For a curve $[C\cup _y D]\in \Delta_2\subset
\mm_5$ where $g(C)=3$ and $g(D)=2$, $\phi([C\cup_y D])$ is a stable
curve of genus $11$ consisting of two disjoint copies $Y_1$ and
$Y_2$ of $C$
 and two disjoint copies $D_1$ and $D_2$ of $D$, such that
$Y_i\cap D_j=\{y_{ij}\}$ for $i, j=1, 2$. The set
$\{y_{1i}, y_{2i}\}\subset D_i$ consists of $y\in D$ and its
hyperelliptic conjugate for each $i=1, 2$. The set $\{y_{1i},
y_{2i}\}\subset Y_i$ for $i=1, 2$, consists of the pairs of points
lying on the tangent line to the smooth plane quartic model of $C$
which passes through the point $y$.
\end{proposition}
\begin{proof} The only case which requires explanation is that when $[C\cup_y D]\in \Delta_2$, when
$\phi([C\cup _y D])$ is the stable reduction of the variety $\overline{G}^1_4(C\cup _y D)$ of limit $\mathfrak g^1_4$'s on $C\cup_y D$. Components of $\overline{G}^1_4(C\cup_y D)$ are indexed by numerical possibilities
 for the ramification
sequences of a limit linear series $l$ such that $\rho(l_C, y)+\rho(l_D, y)=1$ and $\rho(l_C, y), \rho(l_D, y)\geq 0$. When $\rho(l_C, y)=1$ and $\rho(l_D, y)=0$,  we have two numerical possibilities:
\newline
\noindent
\textbf{ (1)}  $a^{l_D}(y)=(1, 4)$, hence $l_D=l_D^1:=y+|\OO_D(3y)|$ and $a^{l_C}(y)\geq (0, 3)$.
Then the curve
$Y_1:=\{l\in G^1_4(C): a^l(y)\geq (0, 3)\}\times \{l_D^1\}$ is an irreducible component of $\phi([C\cup _y D])$.
\newline
\noindent
\textbf{(2)} $a^{l_D}(y)=(2, 3)$, hence $l_D=l_D^2:=2y+\mathfrak g^1_2\in G^1_4(D)$ and $a^{l_C}(y)\geq (1, 2)$.
Then $Y_2:=\{l\in G^1_4(C): a^l(y)\geq (1, 2)\}\times \{l_D^2\}$ is another irreducible component of $\phi([C\cup _y D])$.

Before  we deal with the  remaining case
 when $\rho(l_C, y)=0$ and $\rho(l_D, y)=1$, we note that for a general $[C, y]\in \cM_{3, 1}$, there are two linear
 series $l_C^1, l_C^2 \in G^1_4(C)$ such that $a^{l_C^i}(y)\geq (1, 3)$. If $\rho(l_C, y)=0$, then necessarily
  $a^{l_C}(y)=(1, 3)$, hence $l_C\in \{l_C^1, l_C^2\}$.

We introduce the curves $D_1:=\{l_C^1\}\times \{l\in G^1_4(D): a^l(y)\geq (1, 3)\}$ and  $D_2:=\{l_C^2\}\times  \{l\in G^1_4(D): a^l(y)\geq (1, 3)\}$
which are the remaining two irreducible components of $\phi([C\cup _y D])$.
We single out the points $y_{11}=(l_C^1, l_D^1)\in Y_1\cap D_1,\  y_{12}=(l_C^2, l_D^1)\in Y_1\cap D_2, \ y_{21}=(l_C^1, l_D^2)\in Y_2\cap D_1$
 and $y_{22}=(l_C^2, l_D^2)\in Y_2\cap D_2$ and then $\phi([C\cup_y D])$ is the stable curve of genus $11$ having
  irreducible
 components $Y_1, Y_2, D_1$ and $D_2$ meeting at the points $y_{11}, y_{12}, y_{21}$ and $y_{22}$.
\end{proof}
Proposition \ref{quintic} coupled with Theorem
\ref{divisorspullback} allows us to completely describe the
pull-back map of divisor classes
$\nu^*:\mathrm{Pic}(\overline{\cR}_6)\rightarrow
\mathrm{Pic}(\mm_5)$.

\begin{proposition}\label{prym6} For $\nu: \mm_5 --> \overline{\cR}_6$
given by $[C]\mapsto [W^1_4(C)/ \Gamma]$,  we have the formulas:
$$\nu^*(\lambda)=34\lambda-4\delta_0-15 \delta_1 -(?)\delta_2 ,\mbox{ }
\nu^*(\delta_0^r)=[\gp_{5, 4}^{1, 0}]= 4(33 \lambda-4 \delta_0-15
\delta_1- 21 \delta_2), $$
$$
\nu^*(\delta_0')=\delta_0, \ \mbox{ } \nu^*(\delta_0^{''})=[\mm_{5,
3}^1]=8 \lambda-\delta_0-4 \delta_1-6 \delta_2.
$$
\end{proposition}
\begin{proof} Most of this follows directly by comparing Proposition
\ref{quintic} with the description of the classes $\delta_0^{'},
\delta_0^{''}$ and $\delta_0^r$. Then we use that the generic point
of the Teixidor divisor $\gp_{5, 4}^{1, 0}$ corresponds to a curve
$[C]\in \cM_5$ having precisely one vanishing theta-null (that is,
quadric of rank $3$ containing the canonical image of $C\subset
\PP^4$). In such a case the curve of singular quadrics $[\Gamma
\subset |I_{C/\PP^4}(2)|]\in \mm_6$ is $1$-nodal, the node
corresponding precisely to the vanishing theta-null. This implies
that $\nu^*(\delta_0^r)=[\gp_{5, 4}^{1, 0}]$. Showing that
$\nu^*(\delta_0^{'})=\delta_0$ and $\nu^*(\delta_0^{''})=[\mm_{5,
3}^1]$ proceeds along similar lines. Finally, we write that
$$35 \lambda-4 \delta_0-15 \delta_1- \cdots= \phi^*(\lambda)=\nu^*(\chi^*(\lambda))=
\nu^*\bigl(2\lambda-\frac{1}{4}\ \delta_0^r
\bigr)=2\nu^*(\lambda)-\frac{1}{4}[\gp_{5, 4}^{1, 0}],$$ which
yields the formula for $\nu^*(\lambda)$.
\end{proof}

The main result of \cite{DS} is that the Prym map
$$\mathrm{Prym}:\cR_6\rightarrow \mathcal{A}_5$$ is generically
finite, of degree $27$. We denote by $\mathcal{D}$ the ramification
divisor of $\cR_6\rightarrow \cA_5$ and by $\overline{\mathcal{D}}$
its closure in $\overline{\cR}_6$. It is proved in \cite{B} that the
codifferential  of the Prym map
$$\mathrm{Prym}^*:
T_{\mathrm{Prym} [C, L]}\bigl(\cA_5\bigr)^{\vee} \rightarrow T_{[C,
L]}\bigl(\cR_6\bigr)^{\vee}$$ can be identified with the
multiplication map $$\mbox{Sym}^2 H^0(C, K_C\otimes L)\rightarrow
H^0(C, K_C^{\otimes 2})$$ (Note that $L^{\otimes 2}=\OO_C$). Thus
$[C, L]\in \mathcal{D}$ if and only if $C\stackrel{|K_C\otimes
L|}\longrightarrow \PP^{4}$ lies on a quadric.  An immediate
application of Proposition \ref{prym6} gives the following
characterization of covers of plane quintics which fail the local
Torelli theorem for Pryms:

\begin{theorem}\label{trig}
For the map $\nu: \mm_5 --> \overline{\cR}_6$ given by $[C]\mapsto
[W^1_4(C)/ \Gamma]$, we have the scheme theoretic equality
$\nu^*(\overline{\mathcal{D}})=4\cdot \mm_{5, 3}^1. $ Thus the
abelian variety $\mathrm{Prym}(W^1_4(C)/\Gamma)$ fails the local
Torelli theorem if and only if the curve $[C]\in \cM_5$ is trigonal.
\end{theorem}
\begin{proof} We use  (\ref{prymtauto}) to compute the class of
the compactification $\overline{\mathcal{D}}$ in $\overline{\cR_6}$
of the ramification locus of $\mathrm{Prym}:\cR_6 \rightarrow \cA_5$
(see \cite{F3} for more details and examples). Precisely, there is a
generically non-degenerate morphism between vector bundles of the
same rank $\alpha: \mathrm{Sym}^2(\mathbb E_1)\rightarrow \mathbb
E_2$ over $\overline{\cR}_g$ and
$\overline{\mathcal{D}}=Z_1(\alpha)\cap \cR_g$. From
(\ref{prymtauto}) we find that $$c_1(\mathbb E_2-\mathrm{Sym}^2
(\mathbb
E_1))=7\lambda-\delta_0^{'}-\delta_0^{''}-\frac{3}{2}\delta_0^r-\cdots.$$
By direct computation it follows that
$\nu^*(\overline{\mathcal{D}})\equiv 4\cdot
(8\lambda-\delta_0-a_1\delta_1-a_2\delta_2)$, where $a_1, a_2>1$,
that is, $s(\nu^*(\overline{\mathcal{D}}))=8$. The only irreducible
effective divisor on $\mm_5$ having slope $\leq 8=6+12/(g+1)$ is the
trigonal locus $\mm_{5, 3}^1$, hence we must have the equality of
divisors  $\nu^*(\overline{\mathcal{D}})=4\cdot \mm_{5, 3}^1$.
\end{proof}
\begin{remark} Theorem \ref{trig} is certainly not surprising.
Beauville proves using relatively elementary methods that for any
smooth curve $[C]\in \cM_5-(\cM_{5, 3}^1 \cup \cgp_{5, 4}^{1, 0})$,
the variety $\mathrm{Prym}(W^1_4(C)/\Gamma)$ satisfies local Torelli
(cf. \cite{B}, Proposition 6.4).
\end{remark}

\section{The slope of the moving cone $\mm_g$}

We introduce a fundamental invariant of $\mm_g$ which carries
information about all rational maps from $\mm_g$ to other projective
varieties. If $\mathrm{Mov}(\mm_g)\subset \mathrm{Pic}(\mm_g)\otimes
\mathbb R$ is the cone of moving effective divisors, we define the
\emph{moving slope} of $\mm_g$ by the formula
$$s'(\mm_g):=\mbox{inf}_{D\in \mathrm{Mov}(\mm_g)}\ s(D)\geq s(\mm_g).$$ Any
non-trivial rational map $f: \mm_g-->\PP^N$ provides an upper bound for
$s'(\mm_g)$ because  one has  the obvious inequality $s'(\mm_g)\leq
s\bigl(f^*(\OO_{\PP^N}(1))\bigr)$. This observation is not so useful
for large $g$ when there are very few known examples of rational
maps admitted by $\mm_g$. For low $g$, in the range where $\mm_g$ is
unirational, there are several explicit examples of such maps which
allows us to determine $s'(\mm_g)$. Parts of the next theorem are
certainly known to experts. The slopes $s(\mm_g)$ for $g\leq 11$
have been computed in \cite{FP}, \cite{HMo}, \cite{Ta} and we record
them in the following table for comparison purposes.

\begin{theorem}\label{smallgenus} For integers $3\leq g\leq 11$ we
have the following table concerning the slope and the moving slope
of $\mm_g$ respectively:
\begin{center}
\begin{tabular}{c|cccccccccc}

$g$ & $3$ & $4$& $5$ & $6$ & $7$ & $8$ & $9$ & $10$ & $11$ \\

\hline

$s(\mm_g)$ & $9$ & $\frac{17}{2}$ & $8$ & $\frac{47}{6}$  &
$\frac{15}{2} $ & $\frac{22}{3}$ & $\frac{36}{5}$ & $7$ & $7$ \\

\hline

 $s'(\mm_g)$ & $\frac{28}{3}$ & $[\frac{17}{2}, \frac{44}{5}]$ & $[\frac{41}{5},\frac{33}{4}]$ & $[\frac{47}{6}, \frac{65}{8}]$ & $[\frac{53}{7},
\frac{201}{26}]$ & $[\frac{59}{8}, \frac{149}{20}]$ & $(\frac{36}{5}, \frac{95}{13}]$
& $[\frac{78}{11}, \frac{36}{5}]$ & $7$ \\
\end{tabular}
\end{center}
\end{theorem}

In the proof of Theorem \ref{smallgenus} we use a result, of independent interest, concerning the slopes of curves in  $\mm_g$ which cover the $k$-gonal loci $\cM_{g, k}^1$ for $k\leq 5$. It is a theorem of Tan that if $D\in \mathrm{Eff}(\mm_g)$ is an effective divisor such that $s(D)<7+6/g$ then $D\supset \mm_{g, 3}^1$ (cf. \cite{T}). It is also well-known that if $s(D)<8+4/g$ then $D\supset \mm_{g, 2}^1$ (use that the family arising from a Lefschetz pencil of curves of type $(2, g+1)$ on $\PP^1\times \PP^1$ is a covering curve for $\mm_{g, 2}^1$). Next we prove a similar result for the locus of $4$ and $5$-gonal curves:

\noindent
\emph{Proof of Theorem \ref{4gonal1}.} We begin by recalling that if
$f:X\rightarrow \PP^1$ is a pencil of semi-stable curves of genus $g$ with
$X$ a smooth surface such that there are no $(-1)$-curves in the fibres of $f$,  if $m:\PP^1\rightarrow \mm_g$ denotes
the corresponding moduli map, then the numerical characters of $f$
are computed as follows:
$$\mathrm{deg } \ m^*(\lambda)=\chi(\OO_X)+g-1 \ \mbox{ and }\  \mathrm{deg }\
 m^*(\delta) =c_2(X)+4(g-1).$$ Of
course, these invariants are related by the Noether formula
$12\chi(\OO_X)=K_X^2+c_2(X)$.

The idea of the proof is to use Beniamino Segre's theorem
\cite{S}: A general $k$-gonal curve $[C]\in \cM_{g, k}^1$ has a plane model $\Gamma\subset \PP^2$ of
degree $n\geq (g+k+2)/2$ having one $(n-k)$-fold point $p$ and $$
\delta={n-1\choose 2}-{n-k\choose 2}-g$$
nodes as
the remaining singularities. The pencil $\mathfrak g^1_k$ on $C$ is
recovered by projecting $\Gamma$ from $p$. We denote by
$S:=\mathrm{Bl}_{\delta+1}(\PP^2)$ the blow-up of the plane at $\delta+1$
general points $p_0, \ldots, p_{\delta}\in \PP^2$ and consider the linear system on $S$
$$|\mathcal{L}|=|n\cdot h-(n-k)\cdot
E_{p_0}-2\cdot \sum_{i=1}^{\delta} E_{p_i}|$$
where $h\in \mathrm{Pic}(S)$ is the class of a line. It is known that $|\mathcal{L}|$ is base point free whenever
$$\mbox{virt-dim}(|\mathcal{L}|)=\frac{n(n+3)}{2}-{n-k+1\choose 2}-3\delta\geq 0$$ (cf. \cite{AC2}). This inequality is compatible with the Segre condition precisely when $k\leq 5$, that is, in this range the nodes and the $(n-k)$-fold point of the Segre plane model $\Gamma $ can be chosen to be general points in $\PP^2$.

A covering curve for $\mm_{g, k}^1$ is obtained by blowing-up the $n^2-(n-k)^2-4\delta$ base points of a Lefschetz pencil in the linear system $|\mathcal{L}|$ (see \cite{AC2}, Theorem 5.3 for the fact that one can recover the general curve $[C]\in \cM_{g, k}^1$ in this way). If $F\subset \mm_{g, k}^1$ denotes the induced family, then we have the formulas
$$F\cdot \lambda=g, \ \ F\cdot \delta_0=\frac{k(k+3)}{2}+7g+(3-n)k-3 \mbox{ } \mbox{ and }\  F\cdot \delta_i=0 \mbox{ for all } i\geq 1$$
(For $3\leq k\leq 5$ one checks that there are no $(-1)$-curves in the fibres of $F$, which is not the case for $k=2$). Choosing $n=[(g+k+3)/2]$ (that is, minimal such that Segre's inequality is satisfied), we find that $F\cdot D<0$ implies the inclusion $\mm_{g, k}^1\subset D$ which finishes the proof. Note that for $k=3$ we find that $F\cdot \delta=7g+6$ (independent of $n$), hence $F\cdot \delta/F\cdot \lambda=7+6/g$ and this gives a different proof of Tan's result \cite{T}.
\hfill $\Box$

\begin{corollary}\label{contraction}
There exists no non-trivial rational map $f:\mm_g-->X$ in the projective category such that the indeterminacy locus of $f$ is contained in $\mm_{g, k-1}^1$ and which contracts the variety $\mm_{g, k}^1 (k=4, 5)$  to a point.
\end{corollary}
\begin{proof} By Theorem \ref{4gonal1} we can find two different covering curves $F$ and $F'$ for $\mm_{g, k}^1$ according to different choice of $n\geq (g+k+2)/2$ such that $F\cdot \delta/F\cdot \lambda\neq F'\cdot \delta/F'\cdot \lambda$.\end{proof}
\begin{remark} This last result is in contrast with the situation in the case of the hyperelliptic locus. For instance, the rational map $f:\mm_3-->\mathcal{Q}_3:=|\OO_{\PP^2}(4)|//SL(3)$ to the GIT quotient of plane quartics blows down $\mm_{3, 2}^1$ to the point corresponding to the double conic. Moreover, we have that $f^*(A)\equiv 28\lambda-3\delta-8\delta_1$, where $A\in \mathrm{Ample}(\mathcal{Q}_3)$.
\end{remark}
\noindent
\emph {Proof of Theorem \ref{smallgenus}.}  \textbf{(i) $g=4$.} The Petri divisor $\gp_{4, 3}^1$ is
the closure in $\mm_4$ of the locus of curves $[C]\in \mm_4$ for
which the canonical model of $C\stackrel{|K_C|}\longrightarrow
\PP^3$ lies on a quadric cone. One knows that $\gp_{4, 3}^1\equiv 34
 \lambda- 4 \delta_0- 14 \delta_1- 18 \delta_2$. By taking a
Lefschetz pencil $R\subset \mm_4$ of curves of type $(3, 3)$ on a
smooth quadric in $\PP^3$, we find that $$R\cdot \lambda =4, R\cdot
\delta_0=34$$ which implies that $s(\mm_4)=34/4.$ If $R$ is chosen
generically then $R\cap \gp_{4, 3}^1=\emptyset$. Next we construct a covering curve $F\subset \gp_{4, 3}^1$ for  the Gieseker-Petri divisor. We take the
Hirzebruch surface $\mathbb F_2$ viewed as the blow-up of the cone
$\Lambda\subset \PP^3$ over a conic. We denote as usual,
$\mbox{Pic}(\mathbb F_2)=\mathbb Z\cdot [C_0]\oplus \mathbb Z \cdot
f$, where $f^2=0$, $C_0^2=-2$ and $C_0\cdot f=1$, and $\mathbb
F_2\stackrel{|C_0+2f|}\longrightarrow \PP^3$. Then we consider a
Lefschetz pencil in the linear system $|3C_0+6f|$ corresponding to
intersections of $\Lambda$ with a pencil of cubic surfaces. We
blow-up $\mathbb F_2$ in $18=(3C_0+6f)^2$ base points and denote by
$f:X=\mathrm{Bl}_{18}(\mathbb F_2)\rightarrow \PP^1$ the resulting
family of \emph{semistable} curves. Note that $f$ has precisely one
fibre of the form $C_0+D$ with $D\in |2C_0+6f|$, where $C_0\cdot
D=2$. By blowing-down the $(-2)$-curve $C_0$ we obtain a map
$\nu:X\rightarrow X'$ and a family of \emph{stable} genus $4$ curves
$f':X'\rightarrow \PP^1$, where $X'$ has one surface double point
and $f=f'\circ \nu$.  If $F\subset \gp_{4, 3}^1$ is the curve in the
moduli space induced by $f'$, then $F$ is a covering curve for
$\gp_{4, 3}^1$. Since $\omega_f=\nu^*(\omega_{f'})$, the $\lambda$-
degree of $F$ can be computed directly on $X$, that is,  $F\cdot
\lambda=\chi(\OO_X)+g-1=4$. Then, we can write $F\cdot
\delta=\mathrm{deg} \nu_*([Z])$, where $Z\subset X$ is the $0$-cycle
of nodes in the fibres of $f$, hence
$$F\cdot \delta=12\chi(\OO_X)-K_X^2+4(g-1)=34.$$
Since $F$ and $R$ have the same numerical invariants, it follows that
there is no rational contraction $\mm_4-->X$ having indeterminacy locus contained in $\mm_{4, 2}^1$, which blows the divisor $\gp_{4, 3}^1$
down to a point. The upper bound on $s'(\mm_4)$ is obtained by considering
the irreducible divisor $$\mathfrak{D}_4:=\{[C]\in \cM_4: \exists p\in C \mbox{ with } h^0(C, \OO_C(3p))\geq 2\}$$ introduced by S. Diaz. It is known that $s(\overline{\mathfrak D}_4)=44/5$ (cf. \cite{Di}) , hence $s'(\mm_4)\leq s(\overline{\mathfrak D}_4)$.

\noindent
\textbf{(ii) $g=5$.} We construct a covering curve for
$\gp_{5, 3}^1=\mm_{5, 3}^1$ as follows: On $\mathbb
F_1=\mathrm{Bl}_1(\PP^2)$ we denote by $C_0$ and $f$ respectively,
the generators of the Picard group where $f^2=0, C_0^2=-1, f\cdot
C_0=1$. Then we consider the family of genus $5$ curves $F\subset
\mm_5$ obtained by blowing-up the base points of a  Lefschetz pencil
inside the ample linear system $|3C_0+5f|$ on $\mathbb F_1$. By
direct computation we find $$F\cdot \lambda=5, F\cdot \delta=41,$$
hence $F\cdot \mm_{5, 3}^1=-1$. This implies that $[\mm_{5,
3}^1]\notin \mathrm{Mov}(\mm_5)$ and that $s'(\mm_5)\geq 41/5$. The
upper bound on $s'(\mm_5)$ uses the Teixidor divisor $\gp_{5, 4}^{1, 0}$  which has slope $s(\gp_{5, 4}^{1, 0})=33/4$.

\noindent \textbf{ (iii) $g=6$.}  We use that $s(\gp_{6, 4}^1)=s(\mm_{6})=47/6$ and $s(\gp_{6, 5}^1)=65/8$, hence $s(\gp_{6, 4}^1)\leq s'(\mm_6)\leq s(\gp_{6, 5}^1)$.

\noindent \textbf{ (iv) $g=7$.}  We consider the tetragonal divisor
$\mm_{7, 4}^1\equiv 15
\lambda-2\delta_0-9\delta_1-15\delta_2-18\delta_3$ and we construct
a covering curve for $\mm_{7, 4}^1$ using Theorem \ref{4gonal1}: A general $[C]\in \cM_{7, 4}^1$ has a septic plane model with one triple point and $5$ nodes. A covering curve $F\subset \mm_{7, 4}^1$ is obtained by
blowing up $\PP^2$ at $26=1+5+20$ points, corresponding to the
triple point, the assigned nodes and the unassigned base points of a
Lefschetz pencil in the linear system $$|7\cdot h-3\cdot E_{p_0}-2\cdot \sum_{i=1}^5 E_{p_i}|.$$ We find that $F\cdot \lambda=7, \ F\cdot
\delta=53$, hence $F\cdot \mm_{7, 4}^1<0$. We obtain that $[\mm_{7, 4}^1]\notin \mathrm{Mov}(\mm_7)$  and $$s'(\mm_7)\geq F\cdot
\delta/F\cdot \lambda=\frac{53}{7}.$$

 \noindent \textbf{ (v) $g=8$.} In this
case we consider the Brill-Noether divisor $\mm_{8, 7}^2$
corresponding to septic plane curves with $7$ nodes. To obtain a
covering curve $F\subset \mm_{8, 7}^2$  one has to blow-up $\PP^2$
in the $28=21+7$ base points of a Lefschetz pencil of $7$-nodal
septics. It easily follows that $F\cdot \lambda=8$, $F\cdot
\delta=59$, hence $F \cdot \gp_{8, 7}^2<0$, that is $[\mm_{8,
7}^2]\notin \mathrm{Mov}(\mm_8)$ and $s'(\mm_8)\ge 59/8$. Moreover, $s'(\mm_8)\leq s(\gp_{8, 5}^1)=\frac{149}{20}$ (cf. \cite{F1}).

\noindent \textbf{ (vi) $g=9$.} The smallest known slopes of effective divisors on $\mm_9$ are $s(\mm_{9, 5}^1)=36/5$ and $s(\gp_{9, 8}^2)=95/13$ respectively (cf. \cite{F1}, Theorem 1.5). It follows that a multiple of the linear system $|\gp_{9, 8}^2|$ contains a moving divisor on $\mm_9$.

\noindent \textbf{ (vii) $g=10$.}  We use the results from \cite{FP} and denote by $\overline{\mathcal{K}}_{10}$ the closure of the locus of curves $[C]\in \cM_{10}$ lying on a $K3$ surface, hence $s(\overline{\mathcal{K}}_{10})=7$. If $F\subset \overline{\mathcal{K}}_{10}$ is the $1$-dimensional family obtained from a Lefschetz pencil of curves of genus $10$ lying on a general $K3$ surface, then $F\cdot \delta/F\cdot \lambda=78/11$, hence $s'(\mm_{10})\geq 78/11>s(\overline{\mathcal{K}}_{10})$ and moreover $[\overline{\mathcal{K}}_{10}]\notin \mathrm{Mov}(\mm_{10})$. Since $s(\gp_{10, 6}^1)=36/5$ (cf. \cite{F1}, Proposition 1.6), we obtain the estimate $$\frac{78}{11}\leq s'(\mm_{10})\leq \frac{36}{5}.$$

\noindent \textbf{ (viii) $g=11$.} This is also a consequence of \cite{FP},
Proposition 6.2. If $\overline{\mathcal{F}}_g$ is the Baily-Borel
compactification of the moduli space of polarized $K3$ surfaces of
degree $2g-2$, then there  is a rational map
$f:\mm_{11}-->\overline{\mathcal{F}}_{11}$ given by $f([C])=[S, C]$,
where $S$ is the unique $K3$ surface containing $C$. If $F\subset
\mm_{11}$ is a Lefschetz pencil of curves corresponding to a general
choice of $[S, C]\in \mathcal{F}_{11}$,  then $F\cdot
\lambda=g+1=12$ and $F\cdot \delta=84$. The map $f$ contracts the
pencil $F$, hence for each divisor $A \in
\mbox{Ample}(\overline{\mathcal{F}}_{11})$, we must have that
$s(f^*(A))=7$, that is, $s'(\mm_{11})\leq 7$. Since $F$ is a
covering curve for $\mm_{11}$ one also has that $$s(\mm_{11})\geq \frac{F\cdot
\delta}{ F\cdot \lambda}=7.$$ This gives the estimate $s(\mm_{11})=s'(\mm_{11})=7$.
\hfill $\Box$

\end{document}